\documentclass[11pt,a4paper,reqno]{amsart}
\usepackage[T1]{fontenc}
\usepackage{lmodern}
\usepackage{amsmath,amssymb,mathtools}
\usepackage[margin=28mm]{geometry}
\usepackage{microtype}
\usepackage[numbers,sort&compress]{natbib}
\usepackage[colorlinks=true,linkcolor=blue,citecolor=blue,urlcolor=blue]{hyperref}
\hypersetup{pdftitle={Optimal polynomial decay and nonlinear stabilization for an indirectly damped wave--MGT system},pdfauthor={Tae Gab Ha}}
\numberwithin{equation}{section}
\newtheorem{theorem}{Theorem}[section]
\newtheorem{proposition}[theorem]{Proposition}
\newtheorem{lemma}[theorem]{Lemma}
\newtheorem{corollary}[theorem]{Corollary}
\theoremstyle{remark}
\newtheorem{remark}[theorem]{Remark}
\newcommand{\HH}{\mathcal H}
\newcommand{\GG}{\mathcal G}
\newcommand{\EE}{\mathcal E}
\newcommand{\RR}{\mathbb R}
\newcommand{\CC}{\mathbb C}
\newcommand{\NN}{\mathbb N}
\newcommand{\DD}{\mathcal D}
\newcommand{\ip}[2]{\langle #1,#2\rangle}
\newcommand{\norm}[1]{\lVert #1\rVert}
\DeclareMathOperator{\Rea}{Re}
\DeclareMathOperator{\Ima}{Im}
\DeclareMathOperator{\diag}{diag}

\allowdisplaybreaks[1]

\title[Optimal polynomial decay and nonlinear stabilization for wave--MGT]{Optimal polynomial decay and nonlinear stabilization for an indirectly damped wave--MGT system}
\author{Tae Gab Ha}
\address{Department of Mathematics and Institute of Pure and Applied Mathematics, Jeonbuk National University, Jeonju 54896, Republic of Korea}
\email{tgha@jbnu.ac.kr}
\date{}

\begin{document}
\begin{abstract}
We study a conservative wave equation coupled to a dissipative Moore--Gibson--Thompson equation by a zero-order interaction that is symmetric in augmented variables. A general modal proposition relates the exact resolvent and semigroup asymptotic constants through an asymptotically unitary diagonalization. Its application to this system yields the optimal fourth-order resolvent estimate and the exact linear decay constants. These estimates lead to nonlinear stabilization from small regular initial data. On smooth bounded domains, a weighted composition estimate that accommodates a singular third derivative of the source at zero gives the optimal state rate $t^{-3/4}$ and energy rate $t^{-3/2}$ for power and logarithmic sources. In dimension three, the result includes $f(u)=|u|^q u$ with $3/2<q<2$ and $f(u)=|u|^{\gamma-2}u\log|u|$ with $7/2<\gamma<4$, within the source ranges of the companion general-source theory. Dimension-dependent extensions hold through dimension seven in the same data space, and in arbitrary dimension with higher regularity and a higher order of vanishing at zero. A lower-regularity graph-norm theorem and a conditional transfer theorem for general sources are also proved. Each nonlinear small-data decay exponent is optimal among uniform estimates with the stated dependence on the initial-data norm.
\end{abstract}
\maketitle

\section{Introduction}

Let $\Omega\subset\RR^n$, $n\ge1$, be a bounded domain with $C^2$ boundary. We first consider the linear system
\begin{equation}\label{eq:pde}
\begin{cases}
u_{tt}-\Delta u+\alpha(v+\tau v_t)=0,&\text{in }\Omega\times(0,\infty),\\
\tau v_{ttt}+v_{tt}-\Delta v-b\Delta v_t+\alpha u=0,&\text{in }\Omega\times(0,\infty),\\
u=v=0,&\text{on }\partial\Omega\times(0,\infty).
\end{cases}
\end{equation}
Throughout the paper,
\begin{equation}\label{eq:parameters}
b>\tau>0,\qquad \alpha\in\RR,\qquad 0<|\alpha|<\lambda_1,
\end{equation}
where $\lambda_1$ is the first Dirichlet eigenvalue of $-\Delta$. The wave component has no direct damping. Dissipation is produced by the MGT component and transferred to the wave equation through the terms $\alpha(v+\tau v_t)$ and $\alpha u$.
In the nonlinear problem, the wave equation carries a real source term $f(u)$. Section~\ref{sec:nonlinear-power-log} assumes a $C^\infty$ boundary for its higher-order spatial estimates; the other PDE results retain the $C^2$ assumption. We first determine the sharp linear decay, and then prove nonlinear estimates in explicitly specified source classes and data spaces.

The scalar Moore--Gibson--Thompson equation has been studied through energy and semigroup methods; see \cite{KLM,MMT,BucciEller}. Its high-frequency spectrum contains two oscillatory branches and a bounded real branch. Coupling an MGT field to a conservative elastic field also arises in mixture models. Fern\'andez and Quintanilla \cite{FQ} derive such a system and prove exponential stability under an additional coercivity condition on the principal-order coupling. In \eqref{eq:pde}, the interaction is of order zero, and the two high-frequency propagation speeds are $1$ and $\sqrt{b/\tau}$. The dissipation transfer is correspondingly weaker.

Indirect stabilization of weakly coupled equations has a substantial literature. Alabau, Cannarsa, and Komornik \cite{ACK} developed an abstract framework for indirect internal stabilization, while Alabau-Boussouira and L\'eautaud \cite{AL} studied locally coupled wave-type systems. For systems with memory damping, Oliveira, Oquendo, and Buriol \cite[Theorems~3.3 and~4.1]{OOB} obtained rates depending on the order of the memory operator and the propagation speeds, and proved their optimality for exponential memory kernels. In their unequal-speed case, a memory term of full spatial order yields the state decay exponent $1/4$. Thus this exponent already occurs in closely related hereditary systems. Their formulation uses a history variable and spring-type diagonal and off-diagonal interactions. The present analysis concerns the five-component MGT realization of \eqref{eq:pde} and determines its resolvent and semigroup asymptotic constants in the prescribed energy norm.

The nonlinear version of \eqref{eq:pde} with a logarithmic source was studied in \cite{HaJDE}. That work established the augmented energy structure and the high-frequency wave-branch expansion that excludes uniform exponential stability. The companion paper \cite{HaGeneral} treats general focusing sources and proves, among other results, strong stability of the linearized semigroup. These two papers leave open the global imaginary-axis resolvent estimate and the resulting optimal quantitative decay. We retain their system and energy normalization. The wave-branch expansion is recalled and proved below as part of a complete modal analysis.

Set
\begin{equation}\label{eq:augvariables}
D=b-\tau,\qquad w=v+\tau v_t,\qquad z=v_t,\qquad
Y=(u,u_t,w,w_t,z).
\end{equation}
The energy is
\begin{equation}\label{eq:pdeenergy}
\begin{split}
\EE(Y)=\frac12\bigl(&\norm{u_t}_2^2+\norm{w_t}_2^2
+\norm{\nabla u}_2^2+\norm{\nabla w}_2^2
+\tau D\norm{\nabla z}_2^2\bigr)
+\alpha\Rea(u,w)_{L^2}.
\end{split}
\end{equation}
It is positive definite under \eqref{eq:parameters}, and
\begin{equation}\label{eq:introenergy}
\frac{d}{dt}\EE(Y(t))=-D\norm{\nabla z(t)}_2^2.
\end{equation}
Let $\GG$ be the generator in these variables, $T(t)=e^{t\GG}$, and let $\norm{\cdot}_{E}$ denote the norm satisfying $\norm{Y}_{E}^2=2\EE(Y)$. We use \emph{state decay} for decay of $\norm{Y(t)}_E$ and \emph{energy decay} for decay of the corresponding energy functional. The exact linear results proved here include
\begin{equation}\label{eq:introres}
i\RR\subset\rho(\GG),\qquad
\limsup_{|\omega|\to\infty}|\omega|^{-4}
\norm{(i\omega-\GG)^{-1}}_{E}
=\frac{2(b-\tau)}{\alpha^2},
\end{equation}
and, for every integer $m\ge1$,
\begin{equation}\label{eq:introtime}
\lim_{t\to\infty}t^{m/4}\norm{T(t)\GG^{-m}}_{E}
=\left(\frac{m(b-\tau)}{2e\alpha^2}\right)^{m/4}.
\end{equation}
In particular,
\begin{equation}\label{eq:introbound}
\EE(T(t)Y_0)\le C(1+t)^{-1/2}\norm{\GG Y_0}_{E}^2,
\qquad Y_0\in\DD(\GG),
\end{equation}
with optimal exponent. The constants in \eqref{eq:introres}--\eqref{eq:introtime} refer to this specific energy norm; equivalent norms preserve the exponents, but need not preserve the constants.

The connection between these two exact constants is made explicit in Proposition~\ref{prop:modal-principle}. For an orthogonal sum of finite-dimensional blocks, an asymptotically unitary diagonalization reduces both operator norms to scalar estimates along the spectral branches. The damping coefficients and the leading imaginary parts of the slow branches then determine the constants. In the present system, the required diagonalization follows because the leading modal matrix, after energy normalization and rescaling, is skew-Hermitian with simple spectrum. This eigenvector control is essential: spectral locations alone do not give the resolvent upper bound for a nonnormal operator. Thus the linear contribution consists of the exact asymptotics in the prescribed energy norm and the modal argument connecting them. The Borichev--Tomilov theorem \cite{BT} supplies the corresponding general Hilbert-space relation between resolvent growth and decay rates; see also \cite{BD} for the broader theory of nonuniform stability.

The passage from linear to nonlinear stabilization requires additional work because the decay estimate loses regularity. Semilinear systems indirectly damped through memory were studied by Xiao and Liang \cite{XL}; stabilization of semilinear waves from nonuniform linear decay is also addressed by Joly and Laurent \cite{JL}. Their results concern different systems and hypotheses. In the present problem, the nonlinear source acts in the undamped component, and the energy identity alone gives neither a time-integrable source norm nor a uniform higher-order bound. Our small-data arguments establish both bounds simultaneously by a weighted bootstrap.

The nonlinear conclusions use two initial-data scales. In dimensions $n\le3$, small data in $\DD(\GG)$ yield
\[
\sup_{t\ge0}\norm{Y(t)}_{\DD(\GG)}
+\sup_{t\ge0}(1+t)^{1/4}\norm{Y(t)}_E
\le C\norm{Y_0}_{\DD(\GG)}
\]
under the local vanishing conditions of Section~\ref{sec:nonlinear-graph}: order greater than four for $n=1,2$, and greater than seven for $n=3$. The uniform graph bound is part of the conclusion.

Let $A=-\Delta_D$ be the Dirichlet Laplacian. The main power and logarithmic results use the stronger spatial space
\[
\mathcal X_3=\DD(A^2)\times\DD(A^{3/2})
\times\DD(A^2)\times\DD(A^{3/2})\times\DD(A^2)
\]
and give a uniform $\mathcal X_3$ bound, state decay $t^{-3/4}$, and nonlinear energy decay $t^{-3/2}$. In dimension three they apply to the sources
\[
f(u)=|u|^q u,\quad \tfrac32<q<2,
\qquad
f(u)=|u|^{\gamma-2}u\log|u|,\quad \tfrac72<\gamma<4,
\]
within the strictly $L^2$-subcritical source ranges of \cite{HaGeneral}. Here this condition means that $|f'(\xi)|\le C(1+|\xi|^\varrho)$ for some $\varrho>0$, with $\varrho<2/(n-2)$ when $n\ge3$; in dimensions one and two any finite $\varrho$ is allowed. The nonlinear argument uses a weighted $H^3$ composition estimate that permits a singular third derivative at zero. For absolute powers of degree $\mu$, the sufficient condition $\mu>5/2$ comes from classical Besov composition results \cite{Kateb,BourdaudSharp}. We give a self-contained estimate under local derivative bounds, in a form that also accommodates the logarithmic sources considered here. Spectral Galerkin approximation and energy-space convergence then construct the solution without assuming local Lipschitz continuity of the source map in $\mathcal X_3$.

The same spatial data space allows dimensions $n\le7$, with the dimension-dependent source conditions stated in Section~\ref{sec:nonlinear-power-log}. For every dimension, a further theorem uses $\mathcal X_m$, with integer $m\ge4$, $m\ge n/2$, and sources satisfying a local vanishing condition of order greater than $m$, to obtain state decay $t^{-m/4}$. The precise comparison with the source ranges of \cite{HaGeneral} is made in Section~\ref{sec:nonlinear-power-log}. Each nonlinear exponent is optimal for uniform estimates proportional to the corresponding initial-data norm; the faster rates correspond to stronger initial regularity.

For the general focusing sources of \cite{HaGeneral}, we additionally derive a conditional rate under weighted finite variation of the nonlinear source. This result applies to an existing regular stable trajectory. Remark~\ref{rem:power-log-source-scope} summarizes the scope of the small-data results relative to the companion theory.

We work first with a positive selfadjoint spatial operator with compact inverse. The exact upper limits require no condition on the gaps between its eigenvalues. The full time limit follows when successive eigenvalue ratios tend to one, as they do for the Dirichlet Laplacian by Weyl's law. Section~\ref{sec:abstract} gives the generator and main abstract theorem. Section~\ref{sec:spectral} establishes the modal structure. Section~\ref{sec:resolvent} proves the general modal proposition and applies it to the resolvent; Section~\ref{sec:decay} derives the decay asymptotics and the fixed-datum conclusions. Section~\ref{sec:pde} returns to \eqref{eq:pde} and its original initial data. Section~\ref{sec:nonlinear-graph} treats graph-norm data. Section~\ref{sec:nonlinear-power-log} proves the power and logarithmic results and their higher-dimensional extensions. Section~\ref{sec:nonlinear-variation} gives the weighted-variation theorem and its nonlinear application. Appendix~\ref{app:wave-wave} applies the modal principle to a wave--wave system with unequal propagation speeds.

\section{Energy space, generator, and main theorem}\label{sec:abstract}

Let $X$ be an infinite-dimensional separable complex Hilbert space, with inner product linear in the first variable. Let $A$ be a densely defined positive selfadjoint operator on $X$ with compact inverse. Write
\begin{equation}\label{eq:spatialbasis}
Ae_k=\lambda_k e_k,\qquad
0<\lambda_1\le\lambda_2\le\cdots,\qquad \lambda_k\to\infty,
\end{equation}
where $(e_k)$ is an orthonormal basis and eigenvalues are repeated according to multiplicity. Set $V=\DD(A^{1/2})$, endowed with the norm $\norm{A^{1/2}\cdot}_X$, and let
\[
\HH=V\times X\times V\times X\times V.
\]
For a real-valued problem, we apply the spectral and resolvent analysis to its complexification. All energy estimates restrict to the real space. The real and complex operator norms of $T(t)$ and $T(t)\GG^{-m}$ agree; resolvents at nonreal parameters are understood on the complex space.

For $Y=(u,p_u,w,p_w,z)$, where $p_u$ and $p_w$ denote the velocity components, define
\begin{align}
\norm{Y}_0^2&=\norm{A^{1/2}u}^2+\norm{p_u}^2+
\norm{A^{1/2}w}^2+\norm{p_w}^2+\tau D\norm{A^{1/2}z}^2,\label{eq:norm0}\\
\norm{Y}_E^2&=\norm{Y}_0^2+2\alpha\Rea\ip uw,
\qquad \EE(Y)=\tfrac12\norm{Y}_E^2.\label{eq:normE}
\end{align}
All unmarked component norms are those of $X$. If $\theta=|\alpha|/\lambda_1<1$, then
\begin{equation}\label{eq:equivalence}
(1-\theta)\norm{Y}_0^2\le\norm{Y}_E^2
\le(1+\theta)\norm{Y}_0^2.
\end{equation}
We use $\norm{\cdot}_E$ also for the induced operator norm.

The augmented system is
\begin{equation}\label{eq:abstractsystem}
u_{tt}+Au+\alpha w=0,\qquad
w_{tt}+A(w+Dz)+\alpha u=0,\qquad
\tau z_t+z=w_t.
\end{equation}
The original MGT variable is recovered by $v=w-\tau z$. Indeed, the last equation gives $v_t=w_t-\tau z_t=z$ and hence $w=v+\tau v_t$. Substitution into the middle equation recovers the MGT equation. These identities hold distributionally for energy solutions, and classically whenever the corresponding derivatives exist.
The generator associated with \eqref{eq:abstractsystem} is
\begin{equation}\label{eq:generator}
\GG(u,p_u,w,p_w,z)
=\bigl(p_u,-Au-\alpha w,p_w,-A(w+Dz)-\alpha u,\tau^{-1}(p_w-z)\bigr)
\end{equation}
on
\begin{equation}\label{eq:domain}
\DD(\GG)=\{(u,p_u,w,p_w,z)\in\HH:
p_u,p_w\in V,\ u\in\DD(A),\ w+Dz\in\DD(A)\}.
\end{equation}
The condition in \eqref{eq:domain} is imposed on $w+Dz$, not separately on $w$ and $z$.

\begin{proposition}\label{prop:generation}
Under \eqref{eq:parameters}, $\GG$ generates a contraction $C_0$-semigroup on $(\HH,\norm{\cdot}_E)$. For every $Y_0\in\HH$,
\begin{equation}\label{eq:energyidentity}
\EE(T(t)Y_0)+D\int_0^t\norm{A^{1/2}z(s)}^2\,ds=\EE(Y_0),
\qquad t\ge0.
\end{equation}
\end{proposition}
\begin{proof}
Finite spectral sums are dense in $\HH$ and belong to $\DD(\GG)$. Direct calculation with the inner product associated with \eqref{eq:normE} gives
\begin{equation}\label{eq:dissipation}
\Rea\ip{\GG Y}{Y}_E=-D\norm{A^{1/2}z}^2,
\qquad Y\in\DD(\GG).
\end{equation}
To verify maximality, fix $\xi>0$ and $F=(f_1,\ldots,f_5)\in\HH$. The first, third, and fifth equations of $(\xi-\GG)Y=F$ give
\[
p_u=\xi u-f_1,\qquad p_w=\xi w-f_3,\qquad
z=c_\xi w+h_\xi,
\]
where
\[
c_\xi=\frac{\xi}{1+\tau\xi},\qquad
h_\xi=\frac{\tau f_5-f_3}{1+\tau\xi}\in V.
\]
The remaining equations are the variational system
\begin{equation}\label{eq:laxmilgram}
\begin{split}
&\ip{A^{1/2}u}{A^{1/2}\phi}
+(1+Dc_\xi)\ip{A^{1/2}w}{A^{1/2}\psi}
+\xi^2\bigl(\ip u\phi+\ip w\psi\bigr)
+\alpha\bigl(\ip w\phi+\ip u\psi\bigr)\\
&\hspace{8mm}=\ip{f_2+\xi f_1}\phi+\ip{f_4+\xi f_3}\psi
-D\ip{A^{1/2}h_\xi}{A^{1/2}\psi},
\qquad (\phi,\psi)\in V^2.
\end{split}
\end{equation}
Its form is coercive by $c_\xi\ge0$ and $|\alpha|<\lambda_1$. Lax--Milgram yields a unique $(u,w)\in V^2$. The equations then imply
\[
Au=f_2-\xi p_u-\alpha w\in X,\qquad
A(w+Dz)=f_4-\xi p_w-\alpha u\in X.
\]
The form characterization of $A$ therefore yields $Y\in\DD(\GG)$. The Lumer--Phillips theorem gives generation; see \cite{Pazy}. Identity \eqref{eq:energyidentity} follows first for classical solutions from \eqref{eq:dissipation}, and then for all energy data by approximation, since $z\in C([0,T];V)$ on each bounded time interval.
\end{proof}

For the remainder of the paper set
\begin{equation}\label{eq:a}
a=\frac{\alpha^2}{2D}>0.
\end{equation}

\begin{theorem}\label{thm:main}
Assume \eqref{eq:parameters} and \eqref{eq:spatialbasis}. Then the following assertions hold.
\begin{enumerate}
\item The entire imaginary axis belongs to $\rho(\GG)$, and
\begin{equation}\label{eq:resupper}
\norm{(i\omega-\GG)^{-1}}_E\le C(1+|\omega|)^4,
\qquad \omega\in\RR.
\end{equation}
More precisely,
\begin{equation}\label{eq:exactres}
\limsup_{|\omega|\to\infty}|\omega|^{-4}
\norm{(i\omega-\GG)^{-1}}_E=a^{-1}.
\end{equation}
\item For every integer $m\ge1$,
\begin{equation}\label{eq:mupper}
\norm{T(t)\GG^{-m}}_E\le C_m(1+t)^{-m/4},\qquad t\ge0,
\end{equation}
and
\begin{equation}\label{eq:exactlimsup}
\limsup_{t\to\infty}t^{m/4}\norm{T(t)\GG^{-m}}_E
=\left(\frac{m}{4ae}\right)^{m/4}.
\end{equation}
If additionally
\begin{equation}\label{eq:ratio}
\lambda_{k+1}/\lambda_k\longrightarrow1,
\end{equation}
then the upper limit in \eqref{eq:exactlimsup} is a limit.
\item The semigroup is strongly stable and $\norm{T(t)}_E=1$ for every $t\ge0$. For each fixed $Y_0\in\DD(\GG^m)$,
\begin{equation}\label{eq:individual}
\norm{T(t)Y_0}_E=o(t^{-m/4}).
\end{equation}
\end{enumerate}
\end{theorem}

The positive constant in \eqref{eq:exactlimsup} expresses optimality on the unit ball defined by $\norm{\GG^mY_0}_E\le1$. Since $0\in\rho(\GG)$, this norm is equivalent to the graph norm on $\DD(\GG^m)$. Assertion \eqref{eq:individual} concerns one fixed datum and is consistent with this uniform optimality. No spectral-spacing condition is needed for \eqref{eq:exactres} or \eqref{eq:exactlimsup}.

\section{Modal representation and spectral branches}\label{sec:spectral}

Write $r_k=\sqrt{\lambda_k}$ and $h=\sqrt{D/\tau}$. The map
\begin{equation}\label{eq:modalmap}
\mathcal UY=(X_k)_{k\ge1},\qquad
X_k=(r_ku_k,p_{u,k},r_kw_k,p_{w,k},\sqrt{\tau D}\,r_kz_k)^{\mathsf T},
\end{equation}
is unitary from $(\HH,\norm{\cdot}_0)$ onto $\ell^2(\NN;\CC^5)$. In these coordinates, the generator acts on the $k$th spatial block as $B(r_k)$, where
\begin{equation}\label{eq:modalB}
B(r)=
\begin{pmatrix}
0&r&0&0&0\\
-r&0&-\alpha/r&0&0\\
0&0&0&r&0\\
-\alpha/r&0&-r&0&-hr\\
0&0&0&hr&-1/\tau
\end{pmatrix}.
\end{equation}
The energy on the $k$th block is $|X_k|_{G(r_k)}^2=X_k^*G(r_k)X_k$, with
\begin{equation}\label{eq:metric}
G(r)=I+\frac{\alpha}{r^2}(E_{13}+E_{31}),
\qquad \norm Y_E^2=\sum_k|X_k|_{G(r_k)}^2.
\end{equation}
Here $E_{ij}$ is the elementary matrix. In particular, distinct spatial blocks are also orthogonal with respect to the energy inner product.

\begin{lemma}\label{lem:directsum}
Under $\mathcal U$, the generator is the maximal direct sum of the matrices $B(r_k)$:
\begin{equation}\label{eq:maximaldirectsum}
\mathcal U\DD(\GG)=\{X\in\ell^2:(B(r_k)X_k)_k\in\ell^2\},
\qquad (\mathcal U\GG\mathcal U^{-1}X)_k=B(r_k)X_k.
\end{equation}
Moreover, $(\mathcal UT(t)\mathcal U^{-1}X)_k=e^{tB(r_k)}X_k$.
\end{lemma}
\begin{proof}
The first and third components of $B(r_k)X_k$ are square summable precisely when $p_u,p_w\in V$. The second component imposes $Au+\alpha w\in X$, equivalently $u\in\DD(A)$. The fourth imposes $A(w+Dz)+\alpha u\in X$, equivalently $w+Dz\in\DD(A)$. The fifth is then square summable since $p_w,z\in V$. This proves \eqref{eq:maximaldirectsum}. Each spatial projection preserves the domain and commutes with $\GG$, hence with $T(t)$; the finite-dimensional equation gives the last statement. Alternatively, each modal semigroup is contractive in its energy metric by \eqref{eq:dissipation}, so their direct sum has the same generator.
\end{proof}

\begin{lemma}\label{lem:finite}
For every $k$, all eigenvalues of $B(r_k)$ have strictly negative real part. Consequently, for every fixed $N$ there are constants $C_N,\eta_N>0$ such that
\begin{equation}\label{eq:finiteexp}
\norm{e^{tB(r_k)}}_{G(r_k)}\le C_Ne^{-\eta_Nt},
\qquad 1\le k\le N,\quad t\ge0.
\end{equation}
\end{lemma}
\begin{proof}
Dissipativity places the eigenvalues in the closed left half-plane. If $\GG Y=i\omega Y$ on one spatial block, \eqref{eq:dissipation} gives $z=0$, and the fifth equation gives $p_w=0$. If $\omega\ne0$, the third equation gives $w=0$, the fourth gives $u=0$, and then $p_u=0$. If $\omega=0$, one has $p_u=p_w=z=0$ and
\[
\lambda_ku+\alpha w=0,\qquad \lambda_kw+\alpha u=0.
\]
Since $\lambda_k^2-\alpha^2>0$, $u=w=0$. There are no imaginary eigenvalues. Every finite-dimensional block is therefore exponentially stable, including the case of nontrivial Jordan blocks. Taking a maximum over finitely many blocks proves \eqref{eq:finiteexp}.
\end{proof}

\begin{proposition}\label{prop:spectral}
Let $c_{\mathrm M}=\sqrt{b/\tau}>1$. For sufficiently large real $r$, the five eigenvalues of $B(r)$ are distinct and can be labeled so that, for $\sigma\in\{-1,1\}$,
\begin{align}
s_{\mathrm w,\sigma}(r)
&=\sigma ir-\sigma i a\tau r^{-3}-ar^{-4}+O(r^{-7}),\label{eq:waveexp}\\
s_{\mathrm m,\sigma}(r)
&=\sigma i c_{\mathrm M}r-\frac{D}{2\tau b}+O(r^{-1}),\label{eq:mgtexp}\\
s_0(r)&=-\frac1b-\frac{D}{b^4}r^{-2}+O(r^{-4}).\label{eq:realexp}
\end{align}
There is an eigenvector matrix $S(r)$, converging to a unitary matrix as $r\to\infty$, such that
\begin{equation}\label{eq:diagonal}
B(r)=S(r)\diag(s_{\mathrm w,+},s_{\mathrm w,-},s_{\mathrm m,+},s_{\mathrm m,-},s_0)S(r)^{-1}.
\end{equation}
In the energy metric, its similarity condition number
\begin{equation}\label{eq:kappa}
\kappa(r):=\norm{G(r)^{1/2}S(r)}_2
\norm{S(r)^{-1}G(r)^{-1/2}}_2
\end{equation}
satisfies $\kappa(r)\to1$.
\end{proposition}
\begin{proof}
The modal characteristic equation in the original variables is
\begin{equation}\label{eq:charpoly}
P_r(s)=(s^2+r^2)(\tau s^3+s^2+br^2s+r^2)-\alpha^2(1+\tau s)=0.
\end{equation}
Indeed, the modal ansatz $u(t)=Ue^{st}$, $v(t)=Ve^{st}$ gives the coefficient matrix
\[
\begin{pmatrix}s^2+r^2&\alpha(1+\tau s)\\
\alpha&\tau s^3+s^2+br^2s+r^2\end{pmatrix}.
\]
The determinant of $sI-B(r)$ is $\tau^{-1}P_r(s)$.

Put $\varepsilon=r^{-1}$ and $s=r\widehat{s}$. Dividing \eqref{eq:charpoly} by $r^5$ yields
\begin{equation}\label{eq:F}
F(\widehat{s},\varepsilon)
=\widehat{s}(\widehat{s}^2+1)(\tau\widehat{s}^2+b)
+\varepsilon(\widehat{s}^2+1)^2
-\alpha^2\tau\widehat{s}\varepsilon^4-\alpha^2\varepsilon^5=0.
\end{equation}
At $\varepsilon=0$, its roots are $\pm i,\pm i c_{\mathrm M},0$, all simple. The analytic implicit function theorem gives five analytic roots $\widehat{s}_j(\varepsilon)$ near zero, accounting for the full degree-five polynomial.

For the root through $\widehat{s}_0=\sigma i$, write $\widehat{s}=\widehat{s}_0+\delta$. Since
\[
F_{\widehat{s}}(\sigma i,0)=-2D,\qquad
F(\sigma i,\varepsilon)=-\sigma i\alpha^2\tau\varepsilon^4-\alpha^2\varepsilon^5,
\]
one first obtains $\delta=O(\varepsilon^4)$. Expanding \eqref{eq:F} at $\sigma i$ then gives
\[
0=-2D\delta-\sigma i\alpha^2\tau\varepsilon^4-\alpha^2\varepsilon^5
+O(\delta^2)+O(\varepsilon\delta^2)+O(\varepsilon^4\delta).
\]
Thus
\[
\delta=-\sigma i\frac{\alpha^2\tau}{2D}\varepsilon^4
-\frac{\alpha^2}{2D}\varepsilon^5+O(\varepsilon^8),
\]
which proves \eqref{eq:waveexp}. This recovers the branch used in \cite{HaJDE,HaGeneral} without relying on a formal root approximation.

At $\widehat{s}_0=\sigma i c_{\mathrm M}$,
\[
F_{\widehat{s}}(\sigma i c_{\mathrm M},0)=\frac{2bD}{\tau},\qquad
F_\varepsilon(\sigma i c_{\mathrm M},0)=\frac{D^2}{\tau^2}.
\]
Hence $\widehat{s}_j'(0)=-D/(2\tau b)$, proving \eqref{eq:mgtexp}. For the bounded real root use $\zeta=r^{-2}$ and divide \eqref{eq:charpoly} by $r^4$:
\[
H(s,\zeta)=(1+\zeta s^2)\{bs+1+\zeta(\tau s^3+s^2)\}
-\alpha^2\zeta^2(1+\tau s).
\]
At $(-1/b,0)$ one has $H_s=b$ and $H_\zeta=D/b^3$. The implicit function theorem gives the real analytic root in \eqref{eq:realexp}. Since $s_0/r\to0$, it is the fifth root already constructed.

It remains to control the eigenvectors. One has
\[
C(\varepsilon):=r^{-1}B(r)=B_0+\varepsilon B_1+\varepsilon^2B_2,
\]
where $B_0$ is skew-Hermitian, $B_1=\diag(0,0,0,0,-1/\tau)$, and the only nonzero entries of $B_2$ are $(B_2)_{23}=(B_2)_{41}=-\alpha$. An orthonormal eigenbasis of $B_0$ is
\begin{align*}
v_{\mathrm w,\sigma}&=2^{-1/2}(1,\sigma i,0,0,0)^{\mathsf T},\\
v_{\mathrm m,\sigma}&=2^{-1/2}(0,0,-\sigma i/c_{\mathrm M},1,-\sigma ih/c_{\mathrm M})^{\mathsf T},\\
v_0&=(0,0,-h/c_{\mathrm M},0,1/c_{\mathrm M})^{\mathsf T}.
\end{align*}
For each branch define the spectral projection by the finite product
\begin{equation}\label{eq:projector}
\Pi_j(\varepsilon)=\prod_{\ell\ne j}
\frac{C(\varepsilon)-\widehat{s}_\ell(\varepsilon)I}
{\widehat{s}_j(\varepsilon)-\widehat{s}_\ell(\varepsilon)}.
\end{equation}
The denominators are nonzero for small $\varepsilon$. The vectors
$v_j(\varepsilon)=\Pi_j(\varepsilon)v_j$ are analytic eigenvectors and converge to $v_j$. Taking these vectors as columns gives $S(r)\to S_0$, where $S_0$ is unitary, and proves \eqref{eq:diagonal}. Since $G(r)\to I$, the matrix $G(r)^{1/2}S(r)$ also converges to $S_0$. This proves $\kappa(r)\to1$.
\end{proof}

\begin{corollary}\label{cor:decomposition}
On the complex energy space, the high-frequency eigenvectors together with finitely many generalized eigenvectors form a Riesz basis. There is a bounded projection $P_{\mathrm w}$ commuting with $T(t)$ such that
\begin{equation}\label{eq:fast}
\norm{T(t)(I-P_{\mathrm w})}_E\le Me^{-\eta t},\qquad t\ge0.
\end{equation}
The range of $P_{\mathrm w}$ is the closed span of the eigenvectors associated with the two wave branches in all sufficiently high spatial modes.
\end{corollary}
\begin{proof}
The direct sum of $S(r_k)$ on sufficiently high modes is bounded and boundedly invertible by Proposition~\ref{prop:spectral}. On the remaining finite-dimensional space choose a Jordan basis. Their union is the image of an orthonormal basis under a bounded invertible operator in the product norm, hence a Riesz basis. By \eqref{eq:equivalence}, it is also a Riesz basis in the energy norm. Define $P_{\mathrm w}$ by the first two diagonal coordinates in the high blocks and by zero on the finite complement. Equations \eqref{eq:mgtexp}--\eqref{eq:realexp} and Lemma~\ref{lem:finite} give \eqref{eq:fast}.
\end{proof}

\begin{remark}\label{rem:noncompact}
The compactness of $A^{-1}$ does not imply compactness of the generator resolvent. By \eqref{eq:realexp}, normalized real-branch eigenvectors $Z_k$ satisfy
\[
\GG Z_k=s_0(r_k)Z_k,\qquad s_0(r_k)\to-1/b.
\]
They are orthogonal across spatial modes and bounded in the graph norm of $\GG$. Thus $\DD(\GG)$ does not embed compactly into $\HH$, and $\GG$ has no compact resolvent. The full imaginary-axis resolvent inclusion must therefore be established by uniform estimates, rather than inferred from the absence of imaginary eigenvalues.
\end{remark}

\section{A modal principle and the sharp resolvent estimate}\label{sec:resolvent}

We first isolate the argument that converts modal spectral information into exact operator-norm asymptotics. All norms in the following proposition are taken in the specified Hilbert spaces; the diagonal coordinates carry the Euclidean norm.

\begin{proposition}[Exact constants for asymptotically unitary modal systems]\label{prop:modal-principle}
Let $\mathfrak H=\bigoplus_{k\ge1}\mathfrak H_k$ be an orthogonal sum of complex Hilbert spaces of a fixed finite dimension $N$, and let $0<r_k\le r_{k+1}$ with $r_k\to\infty$. Suppose that every eigenvalue of $B_k:\mathfrak H_k\to\mathfrak H_k$ has negative real part and that, for all sufficiently large $k$, there is an isomorphism $S_k:\CC^N\to\mathfrak H_k$ such that
\begin{equation}\label{eq:general-diagonal}
B_k=S_k\diag(s_{1,k},\ldots,s_{N,k})S_k^{-1},
\qquad \kappa_k:=\norm{S_k}\norm{S_k^{-1}}\longrightarrow1.
\end{equation}
Assume that there are a nonempty set $J\subset\{1,\ldots,N\}$, constants $\beta,\delta>0$, and $a_j>0$, $c_j\in\RR\setminus\{0\}$ for $j\in J$ such that, for high modes,
\begin{align}
\Rea s_{j,k}&=-a_jr_k^{-\beta}+o(r_k^{-\beta}),&
\Ima s_{j,k}&=c_jr_k+o(r_k),&&j\in J,\label{eq:general-slow}\\
\Rea s_{j,k}&\le-\delta,&&&&j\notin J.\label{eq:general-fast}
\end{align}
Let $\mathcal L$ be the maximal direct sum of the $B_k$, namely
\[
\DD(\mathcal L)=\{x\in\mathfrak H:(B_kx_k)_k\in\mathfrak H\},
\qquad \mathcal Lx=(B_kx_k)_k,
\]
and put
\begin{equation}\label{eq:general-constant}
L_0=\max_{j\in J}\frac{1}{a_j|c_j|^\beta}.
\end{equation}
Then $\mathcal L$ generates a bounded $C_0$-semigroup $\mathcal T(t)$,
$i\RR\subset\rho(\mathcal L)$, and
\begin{equation}\label{eq:general-resolvent}
\limsup_{|\omega|\to\infty}|\omega|^{-\beta}
\norm{(i\omega-\mathcal L)^{-1}}=L_0.
\end{equation}
For every integer $m\ge1$,
\begin{equation}\label{eq:general-time}
\limsup_{t\to\infty}t^{m/\beta}\norm{\mathcal T(t)\mathcal L^{-m}}
=\left(\frac{mL_0}{e\beta}\right)^{m/\beta}.
\end{equation}
In particular, the resolvent is $O((1+|\omega|)^\beta)$ and
$\norm{\mathcal T(t)\mathcal L^{-m}}=O((1+t)^{-m/\beta})$.
If $r_{k+1}/r_k\to1$, the upper limit in \eqref{eq:general-time} is a limit.
\end{proposition}
\begin{proof}
For high modes, \eqref{eq:general-diagonal} gives
$\norm{e^{tB_k}}\le\kappa_k$ for $t\ge0$. The remaining finitely many blocks are exponentially stable. Hence the direct sum of the modal semigroups is uniformly bounded, and its strong continuity follows first on finite sums and then by density. Its generator is $\mathcal L$: the coordinate equations give one domain inclusion, and for $x\in\DD(\mathcal L)$ the identity
$\mathcal T(t)x-x=\int_0^t\mathcal T(s)\mathcal Lx\,ds$
gives the other.

For high modes, similarity and normalized eigenvectors give
\begin{equation}\label{eq:general-modal-resolvent}
\max_j|i\omega-s_{j,k}|^{-1}
\le\norm{(i\omega-B_k)^{-1}}
\le\kappa_k\max_j|i\omega-s_{j,k}|^{-1}.
\end{equation}
For fixed $\omega$, the slow-branch terms tend to zero, while the fast-branch terms are at most $\delta^{-1}$. Thus the modal resolvents are uniformly bounded. Their direct sum maps into $\DD(\mathcal L)$, since
$B_k(i\omega-B_k)^{-1}x_k=i\omega(i\omega-B_k)^{-1}x_k-x_k$.
It follows that $i\omega\in\rho(\mathcal L)$ and
\begin{equation}\label{eq:general-resolvent-sup}
\norm{(i\omega-\mathcal L)^{-1}}
=\sup_k\norm{(i\omega-B_k)^{-1}}.
\end{equation}

To obtain the exact upper bound, fix $\vartheta\in(0,1/2)$ and a sufficiently small $\epsilon>0$. For all sufficiently high modes, $\kappa_k\le1+\epsilon$ and
$-\Rea s_{j,k}\ge(a_j-\epsilon)r_k^{-\beta}$ for $j\in J$.
If $(1-\vartheta)|\omega|\le |c_j|r_k\le(1+\vartheta)|\omega|$, then
\[
|\omega|^{-\beta}|i\omega-s_{j,k}|^{-1}
\le\frac{(1+\vartheta)^\beta}{(a_j-\epsilon)|c_j|^\beta}.
\]
Outside this interval, \eqref{eq:general-slow} gives an imaginary-part distance bounded below by a positive multiple of $|\omega|+r_k$, after increasing the modal threshold depending on $\vartheta$. Thus the inverse distance there is $O_\vartheta(|\omega|^{-1})$. The fast-branch contributions are bounded, and each of the finitely many remaining modal resolvents is $O(|\omega|^{-1})$. Equations \eqref{eq:general-modal-resolvent}--\eqref{eq:general-resolvent-sup}, followed by $\epsilon\downarrow0$ and $\vartheta\downarrow0$, give the upper bound $L_0$ in \eqref{eq:general-resolvent}. For the lower bound, take $\omega_{j,k}=\Ima s_{j,k}$ and use an eigenvector to obtain
\[
|\omega_{j,k}|^{-\beta}\norm{(i\omega_{j,k}-\mathcal L)^{-1}}
\ge\frac{|\omega_{j,k}|^{-\beta}}{-\Rea s_{j,k}}
\longrightarrow\frac{1}{a_j|c_j|^\beta}.
\]
Maximizing over $j\in J$ proves \eqref{eq:general-resolvent}.

Since $0\in\rho(\mathcal L)$, similarity gives, for all sufficiently large $k$,
\begin{equation}\label{eq:general-modal-time}
\max_j e^{t\Rea s_{j,k}}|s_{j,k}|^{-m}
\le\norm{e^{tB_k}B_k^{-m}}
\le\kappa_k\max_j e^{t\Rea s_{j,k}}|s_{j,k}|^{-m}.
\end{equation}
The orthogonal direct-sum representation gives
\begin{equation}\label{eq:general-time-sup}
\norm{\mathcal T(t)\mathcal L^{-m}}
=\sup_k\norm{e^{tB_k}B_k^{-m}}.
\end{equation}
The fast branches contribute at most $\delta^{-m}e^{-\delta t}$ before multiplication by $\kappa_k$, and the finite complement also decays exponentially. On each slow branch, $|s_{j,k}|/r_k\to|c_j|$. The elementary identity
\begin{equation}\label{eq:general-maximum}
\sup_{r>0}r^{-m}e^{-dt/r^\beta}
=\left(\frac{m}{e\beta d}\right)^{m/\beta}t^{-m/\beta},
\qquad d,t>0,
\end{equation}
therefore gives, for every sufficiently small $\epsilon>0$,
\[
\limsup_{t\to\infty}t^{m/\beta}\norm{\mathcal T(t)\mathcal L^{-m}}
\le(1+\epsilon)\max_{j\in J}
\frac{1}{(|c_j|-\epsilon)^m}
\left(\frac{m}{e\beta(a_j-\epsilon)}\right)^{m/\beta}.
\]
As $\epsilon\downarrow0$, this is the upper bound in \eqref{eq:general-time}. For the reverse inequality, choose $j$ attaining \eqref{eq:general-constant} and put $t_k=mr_k^\beta/(\beta a_j)$. Then
\[
t_k^{m/\beta}|s_{j,k}|^{-m}e^{t_k\Rea s_{j,k}}
\longrightarrow |c_j|^{-m}\left(\frac{m}{e\beta a_j}\right)^{m/\beta}
=\left(\frac{mL_0}{e\beta}\right)^{m/\beta}.
\]
Equations \eqref{eq:general-modal-time}--\eqref{eq:general-time-sup} prove \eqref{eq:general-time}.

Finally, if $r_{k+1}/r_k\to1$, choose $k=k(t)$ so that
\[
r_{k(t)}\le\left(\frac{\beta a_jt}{m}\right)^{1/\beta}<r_{k(t)+1}.
\]
Then $t/r_{k(t)}^\beta\to m/(\beta a_j)$, and the same eigenvector lower bound holds along all large times. This proves the full limit. The stated global bounds follow from the finite upper limits and boundedness on compact intervals.
\end{proof}

Under the hypotheses of Proposition~\ref{prop:modal-principle}, the two asymptotic constants are therefore related by
\begin{equation}\label{eq:constant-relation}
\limsup_{t\to\infty}t^{m/\beta}\norm{\mathcal T(t)\mathcal L^{-m}}
=\left(\frac{m}{e\beta}
\limsup_{|\omega|\to\infty}|\omega|^{-\beta}
\norm{(i\omega-\mathcal L)^{-1}}\right)^{m/\beta}.
\end{equation}

The diagonalization hypothesis \eqref{eq:general-diagonal} has a useful sufficient condition. If, in orthonormal coordinates on $\mathfrak H_k$, the matrices $r_k^{-1}B_k$ converge to a skew-Hermitian matrix with simple spectrum, their spectral projections converge to its orthogonal rank-one projections. Applying these projections to an orthonormal eigenbasis of the limit produces eigenvector matrices converging to a unitary matrix, hence $\kappa_k\to1$. The damping expansions \eqref{eq:general-slow}--\eqref{eq:general-fast} are separate hypotheses. Proposition~\ref{prop:spectral} verifies both properties for the present system.

\begin{proposition}\label{prop:resupper}
One has $i\RR\subset\rho(\GG)$, the bound \eqref{eq:resupper}, and the exact upper limit \eqref{eq:exactres}.
\end{proposition}
\begin{proof}
Apply Proposition~\ref{prop:modal-principle} to the maximal direct sum in Lemma~\ref{lem:directsum}, with modal metric $G(r_k)$. Lemma~\ref{lem:finite} gives stability of every block, and Proposition~\ref{prop:spectral} gives the required diagonalization. The two wave branches are slow, with
\[
\beta=4,\qquad a_+=a_-=a,\qquad c_+=1,\quad c_-=-1,
\qquad L_0=a^{-1}.
\]
The other branches satisfy \eqref{eq:general-fast} by \eqref{eq:mgtexp}--\eqref{eq:realexp}. Equation~\eqref{eq:general-resolvent} yields the conclusions.
\end{proof}

\begin{remark}
The resolvent upper limit is attained along temporal frequencies aligned with wave-branch eigenvalues. This does not require the spectral-spacing condition \eqref{eq:ratio}, which is used for the full time-domain limit.
\end{remark}

\section{Optimal decay and exact asymptotic constants}\label{sec:decay}

\begin{proposition}\label{prop:timelimsup}
For every integer $m\ge1$, \eqref{eq:mupper} and \eqref{eq:exactlimsup} hold. If \eqref{eq:ratio} holds, then
\begin{equation}\label{eq:exactlimit}
\lim_{t\to\infty}t^{m/4}\norm{T(t)\GG^{-m}}_E
=\left(\frac{m}{4ae}\right)^{m/4}.
\end{equation}
\end{proposition}
\begin{proof}
The modal parameters identified in Proposition~\ref{prop:resupper} give
$\beta=4$ and $L_0=a^{-1}$. Hence \eqref{eq:general-time} gives both \eqref{eq:mupper} and \eqref{eq:exactlimsup}. Since $r_k=\sqrt{\lambda_k}$, condition \eqref{eq:ratio} implies $r_{k+1}/r_k\to1$, so the last assertion of Proposition~\ref{prop:modal-principle} proves \eqref{eq:exactlimit}.
\end{proof}

\begin{proposition}\label{prop:orbits}
The semigroup satisfies assertion~\textup{(3)} of Theorem~\ref{thm:main}.
\end{proposition}
\begin{proof}
Finite spatial sums are dense in $\HH$, and their orbits decay exponentially by Lemma~\ref{lem:finite}. Contractivity and density imply strong stability on $\HH$. On the other hand, normalized wave eigenvectors $Y_k$ give, for every fixed $t\ge0$,
\[
\norm{T(t)}_E\ge e^{t\Rea s_{\mathrm w,+}(r_k)}\longrightarrow1.
\]
Contractivity yields $\norm{T(t)}_E=1$.

Fix $m\ge1$. The operators $t^{m/4}T(t)\GG^{-m}$ are uniformly bounded for $t\ge1$ by \eqref{eq:mupper}, and they converge strongly to zero on finite spatial sums. Density therefore gives
\[
t^{m/4}T(t)\GG^{-m}X\longrightarrow0,\qquad X\in\HH.
\]
Taking $X=\GG^mY_0$ proves \eqref{eq:individual}.
\end{proof}

Propositions~\ref{prop:resupper}, \ref{prop:timelimsup}, and~\ref{prop:orbits} complete the proof of Theorem~\ref{thm:main}. For $m=1$, the implication from \eqref{eq:resupper} to \eqref{eq:mupper} also follows from the Borichev--Tomilov theorem \cite[Theorem~2.4]{BT}, since $T(t)$ is bounded and $i\RR\subset\rho(\GG)$. The direct modal argument supplies, in addition, the exact constants and the spectral-spacing condition for the full time-domain limit.

\begin{remark}\label{rem:graph}
The graph regularity in Theorem~\ref{thm:main} is not a uniform spatial regularity gain in every component. If $\widehat X_{k,j}$ are the high-frequency coordinates after diagonalization, then
\[
Y\in\DD(\GG^m)
\quad\Longleftrightarrow\quad
\sum_{k\text{ high}}\left[
r_k^{2m}\sum_{j\in\{\mathrm w,+;\mathrm w,-;\mathrm m,+;\mathrm m,-\}}
|\widehat X_{k,j}|^2+|\widehat X_{k,0}|^2\right]<\infty,
\]
with the finite-dimensional part unrestricted. This follows from the maximal direct-sum domain and the uniform similarity. The real branch requires no growing spatial weight, but decays exponentially. It is therefore fully compatible with optimal graph-norm decay and with the noncompactness in Remark~\ref{rem:noncompact}.
\end{remark}

\section{Dirichlet problem and initial-data estimates}\label{sec:pde}

Take $X=L^2(\Omega;\CC)$ and $A=-\Delta$ with
\[
\DD(A)=H^2(\Omega)\cap H_0^1(\Omega),\qquad V=H_0^1(\Omega).
\]
The spatial assumptions of Theorem~\ref{thm:main} hold, and Weyl's law gives
\begin{equation}\label{eq:weyl}
\lambda_k\sim (2\pi)^2\left(\frac{k}{\omega_n|\Omega|}\right)^{2/n},
\end{equation}
where $\omega_n$ is the volume of the unit ball in $\RR^n$; see \cite{Davies}. Hence \eqref{eq:ratio} holds.

\begin{corollary}\label{cor:pde}
For \eqref{eq:pde} under \eqref{eq:parameters}, the generator satisfies \eqref{eq:introres}, and the semigroup satisfies \eqref{eq:introtime} for every integer $m\ge1$. Consequently, there are $c_m,C_m>0$ such that
\begin{equation}\label{eq:twosided}
c_m(1+t)^{-m/4}\le\norm{T(t)\GG^{-m}}_E
\le C_m(1+t)^{-m/4},\qquad t\ge0.
\end{equation}
Moreover,
\begin{equation}\label{eq:energyenvelope}
\lim_{t\to\infty}t^{m/2}
\sup_{\substack{Y_0\in\DD(\GG^m)\\\norm{\GG^mY_0}_E\le1}}
\EE(T(t)Y_0)
=\frac12\left(\frac{m(b-\tau)}{2e\alpha^2}\right)^{m/2}.
\end{equation}
\end{corollary}
\begin{proof}
Apply Theorem~\ref{thm:main} and \eqref{eq:weyl}. The positive limit gives \eqref{eq:twosided} for sufficiently large $t$. On a bounded time interval a fixed eigenvector supplies a positive lower bound, while boundedness supplies the upper bound. Finally, $\GG^m:\DD(\GG^m)\to\HH$ is bijective and
\[
\sup_{\norm{\GG^mY_0}_E\le1}\EE(T(t)Y_0)
=\tfrac12\norm{T(t)\GG^{-m}}_E^2,
\]
which proves \eqref{eq:energyenvelope}.
\end{proof}

For the original initial data write
\[
u(0)=u_0,\quad u_t(0)=u_1,\qquad
v(0)=v_0,\quad v_t(0)=v_1,\quad v_{tt}(0)=v_2.
\]
The augmented initial state is
\begin{equation}\label{eq:originalstate}
Y_0=(u_0,u_1,v_0+\tau v_1,v_1+\tau v_2,v_1).
\end{equation}
Energy data are precisely
\[
u_0,v_0,v_1\in V,\qquad u_1,v_2\in X.
\]
For such data, membership in $\DD(\GG)$ is equivalent to
\begin{equation}\label{eq:originaldomain}
u_0\in\DD(A),\qquad u_1,v_2\in V,\qquad
v_0+bv_1\in\DD(A).
\end{equation}
Indeed, $w_0+Dz_0=v_0+bv_1$ and $p_{w,0}=v_1+\tau v_2$.

\begin{corollary}\label{cor:data}
Assume \eqref{eq:originaldomain} in addition to the energy-data conditions. Define
\begin{equation}\label{eq:dataquantity}
\begin{split}
\mathfrak R_0^2={}&\norm{A^{1/2}u_1}^2
+\norm{Au_0+\alpha(v_0+\tau v_1)}^2
+\norm{A^{1/2}(v_1+\tau v_2)}^2\\
&+\norm{A(v_0+bv_1)+\alpha u_0}^2
+\tau D\norm{A^{1/2}v_2}^2.
\end{split}
\end{equation}
Then the solution of \eqref{eq:pde} satisfies
\begin{equation}\label{eq:physicaldecay}
\EE(Y(t))\le C(1+t)^{-1/2}\mathfrak R_0^2,
\qquad t\ge0.
\end{equation}
The exponent $1/2$ cannot be increased in a uniform estimate with this initial-data norm. For each fixed datum satisfying these conditions, $\EE(Y(t))=o(t^{-1/2})$.
\end{corollary}
\begin{proof}
Substitution of \eqref{eq:originalstate} into \eqref{eq:generator} gives
\[
\GG Y_0=\bigl(u_1,-Au_0-\alpha(v_0+\tau v_1),
v_1+\tau v_2,-A(v_0+bv_1)-\alpha u_0,v_2\bigr).
\]
Thus $\mathfrak R_0=\norm{\GG Y_0}_0$, which is equivalent to $\norm{\GG Y_0}_E$ by \eqref{eq:equivalence}. The claims follow from Theorem~\ref{thm:main} and Corollary~\ref{cor:pde}.
\end{proof}

\begin{remark}\label{rem:scope}
All asymptotics are for fixed parameters satisfying \eqref{eq:parameters}. When $\alpha=0$, the wave component is conservative. When $b=\tau$, the energy dissipation vanishes and the two oscillatory speeds coincide, so the simple-root diagonalization used above no longer applies. Thus the estimates depend on the fixed parameters and may degenerate in either limit. The same convention applies to the nonlinear results below.
\end{remark}

\section{Sharp decay for small nonlinear solutions from graph-norm data}
\label{sec:nonlinear-graph}

In this section $\Omega\subset\RR^n$, $n\in\{1,2,3\}$, is a bounded domain
with $C^2$ boundary, and $A=-\Delta$ is its Dirichlet realization. We use
the real energy space and retain $b>\tau>0$ and $0<|\alpha|<\lambda_1$.
Consider the semilinear system
\begin{equation}\label{eq:nonlinear-pde}
\begin{cases}
u_{tt}-\Delta u+\alpha(v+\tau v_t)=f(u),\\
\tau v_{ttt}+v_{tt}-\Delta v-b\Delta v_t+\alpha u=0,
\end{cases}
\end{equation}
with the same homogeneous Dirichlet conditions as before. Write
$\mathcal N(Y)=(0,f(u),0,0,0)$, so that its augmented formulation is
$Y_t=\GG Y+\mathcal N(Y)$.

We impose only a local condition on the source: there exist
$r_0>0$, $C_f>0$, and an exponent $\mu$, used throughout the nonlinear estimates to describe the local vanishing condition, satisfying
\begin{equation}\label{eq:nonlinear-order}
 \mu>4 \quad(n=1,2),\qquad \mu>7 \quad(n=3),
\end{equation}
such that $f$ is real valued and $C^2$ on
an open interval containing $[-r_0,r_0]$, and
\begin{equation}\label{eq:nonlinear-source}
 |f^{(j)}(\xi)|\le C_f|\xi|^{\mu-j},\qquad |\xi|\le r_0,
 \quad j=0,1,2.
\end{equation}
In particular, $f(0)=f'(0)=0$. No condition on the behavior of $f$ at
infinity enters the result below. One may extend $f$ outside this
interval as a $C^2$ function, or multiply such an extension by a smooth
cutoff. The solutions constructed below stay within the interval on
which the original source is prescribed.

Let
\begin{equation}\label{eq:nonlinear-graphnorm}
 \mathcal Z=\DD(\GG),\qquad \norm{Y}_{\mathcal Z}=\norm{\GG Y}_E.
\end{equation}
Since $0\in\rho(\GG)$, this is a complete norm equivalent to the usual
graph norm. The restricted semigroup is a contraction on $\mathcal Z$ because
$\GG T(t)Y=T(t)\GG Y$. Its strong continuity follows from the strong
continuity of $T(t)$ on $\HH$. Theorem~\ref{thm:main} gives
\begin{equation}\label{eq:nonlinear-linear-bounds}
 \norm{T(t)Z}_{\mathcal Z}\le\norm Z_{\mathcal Z},\qquad
 \norm{T(t)Z}_E\le C(1+t)^{-1/4}\norm Z_{\mathcal Z},
 \qquad Z\in \mathcal Z.
\end{equation}

\begin{theorem}\label{thm:nonlinear-small}
Assume \eqref{eq:nonlinear-order}--\eqref{eq:nonlinear-source}. There exist $\varepsilon_*>0$ and
$C>0$ such that every real $Y_0\in \mathcal Z$ with
$\norm{Y_0}_{\mathcal Z}\le\varepsilon_*$ generates a unique global solution of
\eqref{eq:nonlinear-pde} satisfying
\begin{equation}\label{eq:nonlinear-regularity}
 Y\in C([0,\infty);\mathcal Z)\cap C^1([0,\infty);\HH),
 \qquad \sup_{t\ge0}\norm{u(t)}_{L^\infty}\le r_0/2,
\end{equation}
and
\begin{equation}\label{eq:nonlinear-decay}
 \sup_{t\ge0}\norm{Y(t)}_{\mathcal Z}
 +\sup_{t\ge0}(1+t)^{1/4}\norm{Y(t)}_E
 \le C\norm{Y_0}_{\mathcal Z}.
\end{equation}
In addition,
\begin{equation}\label{eq:nonlinear-remainder}
 \norm{Y(t)-T(t)Y_0}_E
 \le C\norm{Y_0}_{\mathcal Z}^\mu(1+t)^{-1/4}.
\end{equation}
Let $F(\xi)=\int_0^\xi f(\zeta)\,d\zeta$ and define
\begin{equation}\label{eq:nonlinear-energy}
 \EE_f(Y)=\EE(Y)-\int_\Omega F(u)\,dx.
\end{equation}
Along these solutions this energy is positive definite, uniformly
equivalent to $\norm Y_E^2$, and satisfies
\begin{equation}\label{eq:nonlinear-energy-decay}
 \EE_f(Y(t))+D\int_0^t\norm{\nabla z(s)}_2^2\,ds
 =\EE_f(Y_0),\qquad
 \EE_f(Y(t))\le C(1+t)^{-1/2}\norm{Y_0}_{\mathcal Z}^2.
\end{equation}
The state exponent $1/4$, and hence the energy exponent $1/2$, is
optimal among uniform estimates with the stated dependence on the initial graph norm
on a neighborhood of zero in $\mathcal Z$.
\end{theorem}

\begin{proof}
We first record the regularity actually supplied by the graph norm.
The second component of $\GG Y$ yields
$Au=-(\GG Y)_2-\alpha w$. Consequently, elliptic regularity and
$\norm Y_E\le\norm{\GG^{-1}}_E\norm Y_{\mathcal Z}$ give
\begin{equation}\label{eq:nonlinear-uH2}
 \norm u_{H^2}\le C\norm Y_{\mathcal Z},\qquad
 \norm u_{H^1}\le C\norm Y_E.
\end{equation}
Only $u$ is used in this estimate. In particular, it does not require
separate $H^2$ estimates for $w$ and $z$, or compactness of the generator
resolvent.

For $\norm u_\infty\le r_0$, the trace condition and $f(0)=0$ imply
$f(u)\in H_0^1(\Omega)$. Thus $\mathcal N(Y)\in \mathcal Z$ and
\[
 \GG\mathcal N(Y)=(f(u),0,0,0,0),\qquad
 \norm{\mathcal N(Y)}_{\mathcal Z}=\norm{\nabla f(u)}_2.
\]
It follows that
\begin{equation}\label{eq:nonlinear-tame}
 \norm{\mathcal N(Y)}_{\mathcal Z}
 \le C\norm u_\infty^{\mu-1}\norm u_{H^1}.
\end{equation}
The chain rule and the bound on $f''$ also give, whenever
$\norm u_\infty,\norm{\widetilde u}_\infty\le r_0$,
\begin{align}\label{eq:nonlinear-lipschitz}
 \norm{\mathcal N(Y)-\mathcal N(\widetilde Y)}_{\mathcal Z}
 &\le C(\norm u_\infty+\norm{\widetilde u}_\infty)^{\mu-1}
       \norm{u-\widetilde u}_{H^1}\notag\\
 &\quad+C(\norm u_\infty+\norm{\widetilde u}_\infty)^{\mu-2}
       \norm{u-\widetilde u}_\infty\norm{\widetilde u}_{H^1}.
\end{align}
Since $H^2(\Omega)\hookrightarrow L^\infty(\Omega)$ for $n\le3$,
this proves local Lipschitz continuity on $\mathcal Z$ in the region under
consideration. Standard local semilinear theory on the Banach space
$\mathcal Z$ therefore gives a unique mild solution satisfying
\begin{equation}\label{eq:nonlinear-duhamel}
 Y(t)=T(t)Y_0+\int_0^tT(t-s)\mathcal N(Y(s))\,ds.
\end{equation}
The $\mathcal Z$-valued continuity of the source and $Y_0\in\DD(\GG)$ imply
$Y\in C^1(\HH)$ on its interval of existence. Alternatively, this
follows by differentiating the integral in \eqref{eq:nonlinear-duhamel}
in $\HH$, using
$\GG\mathcal N(Y)\in C(\HH)$.

In dimension one set $\delta=0$. In dimensions two and three choose, respectively,
\begin{equation}\label{eq:nonlinear-delta}
 0<\delta<\frac{\mu-4}{\mu-1}\quad(n=2),\qquad
 \frac12<\delta<\frac{\mu-4}{\mu-1}\quad(n=3).
\end{equation}
The intervals are nonempty precisely under the strict inequalities in
\eqref{eq:nonlinear-order}.
The Sobolev embedding and interpolation inequalities give
\begin{equation}\label{eq:nonlinear-interpolation}
 \norm u_\infty\le C\norm u_{H^1}^{1-\delta}
                         \norm u_{H^2}^{\delta}.
\end{equation}
For $n=1$ this is just $H^1(\Omega)\hookrightarrow L^\infty(\Omega)$;
for $n=2,3$ use $H^{1+\delta}(\Omega)\hookrightarrow L^\infty(\Omega)$,
since $1+\delta>n/2$.
Put
\begin{equation}\label{eq:nonlinear-exponent}
 \nu_f=1+(\mu-1)(1-\delta)>4,\qquad \eta=\nu_f/4>1.
\end{equation}
For $T$ below the maximal existence time define
\[
 M_T=\sup_{0\le t\le T}\norm{Y(t)}_{\mathcal Z}
       +\sup_{0\le t\le T}(1+t)^{1/4}\norm{Y(t)}_E.
\]
As long as the solution remains in the prescribed source interval,
\eqref{eq:nonlinear-uH2}--\eqref{eq:nonlinear-interpolation} imply
\begin{equation}\label{eq:nonlinear-force-bound}
 \norm{\mathcal N(Y(t))}_{\mathcal Z}
 \le C M_T^\mu(1+t)^{-\eta}.
\end{equation}
The two linear estimates in \eqref{eq:nonlinear-linear-bounds}, applied
to \eqref{eq:nonlinear-duhamel}, now show
\begin{align*}
 \norm{Y(t)}_{\mathcal Z}
 &\le\norm{Y_0}_{\mathcal Z}+C M_T^\mu\int_0^t(1+s)^{-\eta}\,ds,\\
 \norm{Y(t)}_E
 &\le C(1+t)^{-1/4}\norm{Y_0}_{\mathcal Z}
  +C M_T^\mu\int_0^t(1+t-s)^{-1/4}(1+s)^{-\eta}\,ds.
\end{align*}
Since $\eta>1$,
\begin{equation}\label{eq:nonlinear-convolution}
 \int_0^t(1+t-s)^{-1/4}(1+s)^{-\eta}\,ds
 \le C_\eta(1+t)^{-1/4}.
\end{equation}
Indeed, on $[0,t/2]$ the first factor is bounded by
$C(1+t)^{-1/4}$ and the second is integrable. On $[t/2,t]$ the integral
is at most $C(1+t)^{-\eta+3/4}$, which is bounded by the right-hand side
because $\eta>1$.
We have therefore proved
\begin{equation}\label{eq:nonlinear-bootstrap}
 M_T\le C_0\norm{Y_0}_{\mathcal Z}+C_1M_T^\mu.
\end{equation}
Choose $C_0$ also larger than the embedding constant controlling
$M_0/\norm{Y_0}_{\mathcal Z}$. By decreasing $\varepsilon_*$ if necessary, a
continuity argument applied to \eqref{eq:nonlinear-bootstrap} gives
$M_T\le2C_0\norm{Y_0}_{\mathcal Z}$ for every admissible $T$. The same choice can
ensure, by \eqref{eq:nonlinear-uH2}, that
$\norm{u(t)}_\infty\le r_0/2$. Thus the solution cannot leave the
interval on which \eqref{eq:nonlinear-source} holds. Local continuation
in $\mathcal Z$ rules out a finite maximal time: the $\mathcal Z$ norm remains bounded and
the nonlinear map is Lipschitz on the corresponding bounded set after
the harmless extension described above. This proves global existence
and \eqref{eq:nonlinear-decay}. Applying the second linear estimate to
the integral alone, with the established bound on $M_T$, proves
\eqref{eq:nonlinear-remainder}.

The potential satisfies $|F(\xi)|\le C|\xi|^{\mu+1}$. Hence
\[
 \left|\int_\Omega F(u)\,dx\right|
 \le C\norm u_\infty^{\mu-1}\norm u_2^2
 \le C\norm u_\infty^{\mu-1}\norm Y_E^2.
\]
Further reducing $\varepsilon_*$ yields, for example,
\begin{equation}\label{eq:nonlinear-energy-equivalence}
 \tfrac14\norm{Y(t)}_E^2
 \le\EE_f(Y(t))\le\tfrac34\norm{Y(t)}_E^2.
\end{equation}
The energy identity follows by taking the energy inner product of the
equation with $Y$ and using
\[
 \frac{d}{dt}\int_\Omega F(u(t))\,dx=(f(u(t)),u_t(t))_{L^2}.
\]
All terms are defined at the regularity already proved. Integration
gives \eqref{eq:nonlinear-energy-decay}, whose decay estimate follows
from \eqref{eq:nonlinear-decay} and
\eqref{eq:nonlinear-energy-equivalence}.

Finally, let $\mathcal S(t)Y_0$ denote the nonlinear solution. The local $C^2$
assumption on $f$ makes $\mathcal N:\mathcal Z\to \mathcal Z$ continuously
Fr\'echet differentiable on the small neighborhood under consideration,
with $D\mathcal N(0)=0$. Alternatively, the explicit remainder estimate
gives, for every $Z\in \mathcal Z$ and every fixed $t\ge0$,
\begin{equation}\label{eq:nonlinear-tangent}
 \lim_{\varepsilon\to0}
 \varepsilon^{-1}\mathcal S(t)(\varepsilon Z)=T(t)Z
 \quad\hbox{in }\HH.
\end{equation}
Thus the derivative of the flow at zero is the linear semigroup.
If, on any graph-norm neighborhood of zero, there were a uniform bound
\[
 \norm{\mathcal S(t)Y_0}_E\le C(1+t)^{-\beta}\norm{Y_0}_{\mathcal Z},
 \qquad t\ge0,\quad \beta>1/4,
\]
substitution of $Y_0=\varepsilon Z$, division by $|\varepsilon|$, and
\eqref{eq:nonlinear-tangent} would imply the same bound for $T(t)Z$.
This contradicts the sharp linear decay theorem, with precisely the
norm $\norm Z_{\mathcal Z}=\norm{\GG Z}_E$. The same argument excludes an
improvement by any common factor tending to zero multiplying
$(1+t)^{-1/4}$. A faster uniform nonlinear energy estimate would, by
\eqref{eq:nonlinear-energy-equivalence}, imply a faster state estimate.
This proves the stated optimality.
\end{proof}

\begin{corollary}\label{cor:nonlinear-examples}
Theorem~\ref{thm:nonlinear-small} applies to the following sources
and dimensions:
\begin{enumerate}
\item $f(\xi)=|\xi|^q\xi$ with $q>3$, in dimensions $n=1,2$;
\item $f(\xi)=|\xi|^{\gamma-2}\xi\log|\xi|$, $f(0)=0$, with $\gamma>5$,
      in dimensions $n=1,2$;
\item in dimension $n=3$, the source
\begin{equation}\label{eq:nonlinear-flattened-source}
 f(\xi)=\frac{\xi^9}{(1+\xi^2)^{7/2}}.
\end{equation}
\end{enumerate}
\end{corollary}
\begin{proof}
For the power source take $\mu=q+1$. For the logarithmic source choose
$4<\mu<\gamma-1$. For $j=0,1,2$ its derivatives away from zero are
bounded near zero by
$C|\xi|^{\gamma-1-j}(1+|\log|\xi||)$, which is at most
$C_\mu|\xi|^{\mu-j}$. These derivatives extend continuously through zero,
so \eqref{eq:nonlinear-source} holds on a sufficiently small interval.
The third source is smooth and satisfies \eqref{eq:nonlinear-source}
near zero with $\mu=9>7$. It also belongs to the three-dimensional
$L^2$-subcritical focusing class of the general-source companion paper.
Indeed,
\[
 f'(\xi)=\frac{\xi^8(9+2\xi^2)}{(1+\xi^2)^{9/2}},\qquad
 |f'(\xi)|\le C(1+|\xi|),
\]
so its global derivative growth exponent is $1<2=2/(n-2)$.
For $\xi\ne0$,
\[
 \frac{\xi f'(\xi)}{f(\xi)}=2+\frac{7}{1+\xi^2}\ge2.
\]
With $H_3(\xi)=\xi f(\xi)/3-F(\xi)$, this gives, for $\xi>0$,
\[
 \xi H_3'(\xi)=\frac{\xi^2f'(\xi)-2\xi f(\xi)}{3}
         =\frac{7\xi f(\xi)}{3(1+\xi^2)}>0.
\]
Since $H_3$ is even and $H_3(0)=0$, it is nonnegative and radially
nondecreasing. Thus the structural hypotheses hold with $\theta=3$.
Moreover, $F(\xi)>0$ for $\xi\ne0$, so a nonzero smooth compactly supported
test function satisfies the positive-potential hypothesis. This
example distinguishes the high order of vanishing at the origin,
which yields small-data decay, from the subcritical growth at infinity.
\end{proof}

\begin{remark}\label{rem:nonlinear-scope}
The source restrictions ensure integrability of the graph-valued
nonlinear source through $\eta>1$ in
\eqref{eq:nonlinear-exponent}. The endpoints $\mu=4$ for $n=1,2$
and $\mu=7$ for $n=3$ are therefore excluded by this argument.
The examples above illustrate that the local vanishing condition
and growth at infinity play different roles. Section~\ref{sec:nonlinear-power-log}
uses stronger spatial data norms to treat lower-order power and
logarithmic sources. The nonlinear conclusions concern optimal
uniform exponents; the exact constants in Theorem~\ref{thm:main}
refer to the linear semigroup.
\end{remark}

\begin{remark}[Dimension and choice of data space]
The restriction $n\le3$ in this section comes from
$H^2(\Omega)\hookrightarrow L^\infty(\Omega)$, used to control the
local source by the graph norm. This embedding fails for $n\ge4$,
so the proof does not extend by changing the dimension alone.
Section~\ref{sec:nonlinear-power-log} uses stronger spatial norms
and proves higher-dimensional results. The linear estimates in
Sections~\ref{sec:abstract}--\ref{sec:pde} have no such dimension restriction.
\end{remark}

\section{Small-data stabilization for power and logarithmic sources}
\label{sec:nonlinear-power-log}

Throughout this section, $A=-\Delta_D$ on a bounded $C^\infty$ domain
$\Omega\subset\RR^n$, and all functions and phase spaces are real.
The coefficients satisfy $D=b-\tau>0$ and
$0<|\alpha|<\lambda_1$. We consider
\begin{equation}\label{eq:spatial-nonlinear-system}
\begin{cases}
 u_{tt}+Au+\alpha w=f(u),\\
 w_{tt}+A(w+Dz)+\alpha u=0,\\
 \tau z_t+z=w_t.
\end{cases}
\end{equation}
The estimates below use a spatial scale, rather than powers of the
generator. This distinction is necessary because the bounded real
spectral branch prevents an identification of these two scales.

\subsection{Spatial decay and composition estimates}

For $s\ge0$, define
\begin{equation}\label{eq:spatial-scale}
\begin{split}
\mathcal X_s={}&\DD(A^{(s+1)/2})\times\DD(A^{s/2})
\times\DD(A^{(s+1)/2})\times\DD(A^{s/2})
\times\DD(A^{(s+1)/2}),\\
\norm{Y}_{\mathcal X_s}={}&\norm{\mathbf A^{s/2}Y}_E,
\end{split}
\end{equation}
where $\mathbf A^{s/2}$ acts componentwise. Thus
$\mathcal X_0=\HH$, and the equivalent product norm in
\eqref{eq:spatial-scale} contains $s+1$ spatial derivatives of
$u,w,z$ and $s$ derivatives of $u_t,w_t$.

\begin{lemma}[Linear decay in the spatial scale]\label{lem:spatial-decay}
For every $s\ge0$, the linear semigroup is a strongly continuous
contraction semigroup on $\mathcal X_s$, and
\begin{equation}\label{eq:spatial-decay}
 \norm{T(t)Y}_E\le C_s(1+t)^{-s/4}\norm{Y}_{\mathcal X_s},
 \qquad t\ge0,\quad Y\in\mathcal X_s.
\end{equation}
For $s>0$, the uniform exponent $s/4$ is optimal, including with
respect to improvements by a common factor tending to zero.
\end{lemma}
\begin{proof}
The spectral projections of $A$ commute with $T(t)$. Applying the
energy identity to finite spectral sums of $\mathbf A^{s/2}Y$,
and then using density, proves contractivity and strong continuity
on $\mathcal X_s$.
Write $r_k=\sqrt{\lambda_k}$. Uniform modal diagonalization and
the spectral branches established above give
\[
 \norm{T_k(t)}_{E,k}\le C\exp(-ct/r_k^4),\qquad t\ge0,
\]
with constants independent of $k$. Indeed, the wave branches have
real parts at most $-cr_k^{-4}$, and the remaining branches are
uniformly separated from the imaginary axis. The finitely many
remaining modes are exponentially stable; possible Jordan factors
are absorbed by decreasing the exponential constant.
Since the $\mathcal X_s$ norm on the $k$th mode is $r_k^s$
times its energy norm, \eqref{eq:spatial-decay} follows from
\[
 \sup_{r\ge r_1}r^{-s}e^{-ct/r^4}\le C_s(1+t)^{-s/4}.
\]
For optimality, take a wave eigenvector $V_k$ with
$\norm{V_k}_E=1$ and eigenvalue $\zeta_k$ satisfying
$-\Rea\zeta_k\sim a r_k^{-4}$, where $a=\alpha^2/(2D)$.
The vector $Z_k=r_k^{-s}V_k$ has unit $\mathcal X_s$ norm,
and at $t_k=(-\Rea\zeta_k)^{-1}$,
\[
 \norm{T(t_k)Z_k}_E=e^{-1}r_k^{-s}\asymp t_k^{-s/4}.
\]
The real and complex operator norms agree for the complexification
of an operator between real Hilbert spaces, so these eigenvectors
also establish the assertion for the real spaces used here.
\end{proof}

The third-order composition estimate below allows the third derivative
of the source to be singular at the origin. For absolute powers,
Kateb~\cite{Kateb} proved boundedness on Besov spaces; in particular,
$u\mapsto|u|^\mu$ is bounded on $B^s_{2,2}\cap L^\infty$
for $\mu>1$ and $0<s<\mu+1/2$.
See also \cite{BourdaudSharp} for sharpness and \cite{Bourdaud}
for general background. Since $H^3=B^3_{2,2}$, this gives the
familiar sufficient threshold $\mu>5/2$.
We prove a self-contained weighted tame estimate under local
derivative bounds, which applies to both the signed powers and
the logarithmic sources in this section.

\begin{lemma}[A real third-order composition estimate]
\label{lem:real-power-composition}
Let $\mu>5/2$ and $r_0>0$. Suppose that
\[
 f\in C^2((-r_0,r_0))\cap C^3((-r_0,r_0)\setminus\{0\}),
 \qquad f(0)=f'(0)=f''(0)=0,
\]
and that
\begin{equation}\label{eq:spatial-source-assumption}
 |f^{(j)}(\xi)|\le C_f|\xi|^{\mu-j},
 \qquad 0<|\xi|<r_0,\quad j=0,1,2,3.
\end{equation}
On every bounded smooth domain, there exists $r_*>0$ such that
\begin{equation}\label{eq:real-composition-tame}
 \norm{f(u)}_{H^3(\Omega)}
 \le C\norm u_\infty^{\mu-1}\norm u_{H^3(\Omega)}
\end{equation}
whenever $u\in H^3(\Omega)\cap L^\infty(\Omega)$ is real and
$\norm u_\infty\le r_*$. The constants may depend on the domain,
$\mu$, and the source bounds, but not on $u$.
\end{lemma}
\begin{proof}
We first work on $\RR^n$ with a real compactly supported smooth
function $U$, and write $M=\norm U_\infty$. The case $M=0$ is
immediate; otherwise, assume initially that $M<r_0/2$.
For $5/2<\mu<3$, put
\[
 W_\varepsilon(\xi)=(\xi^2+\varepsilon^2)^{\mu-3},
 \qquad
 \Psi_\varepsilon(\xi)=\int_0^\xi W_\varepsilon(\zeta)\,d\zeta.
\]
Since $2\mu-6>-1$, scaling in this integral shows that
\begin{equation}\label{eq:primitive-weight-bound}
 |\Psi_\varepsilon(\xi)|\le C_\mu|\xi|W_\varepsilon(\xi).
\end{equation}
For a coordinate direction $i$, define
\[
 I_{i,\varepsilon}
 =\int_{\RR^n}W_\varepsilon(U)|\partial_iU|^6\,dx.
\]
Integration by parts, followed by Cauchy--Schwarz, gives
\[
 \begin{split}
 I_{i,\varepsilon}
 &=-5\int_{\RR^n}\Psi_\varepsilon(U)
           (\partial_iU)^4\partial_{ii}U\,dx\\
 &\le5I_{i,\varepsilon}^{1/2}
 \left(\int_{\RR^n}
 \frac{|\Psi_\varepsilon(U)|^2}{W_\varepsilon(U)}
 |\partial_iU|^2|\partial_{ii}U|^2\,dx\right)^{1/2}.
 \end{split}
\]
The quotient in the last integral satisfies
\[
 \frac{|\Psi_\varepsilon(\xi)|^2}{W_\varepsilon(\xi)}
 \le C_\mu|\xi|^2(\xi^2+\varepsilon^2)^{\mu-3}
 \le C_\mu|\xi|^{2\mu-4}.
\]
The Gagliardo--Nirenberg derivative interpolation inequalities
\cite[Lecture~II, p.~125]{Nirenberg} imply
\[
 \norm{\partial_iU}_6
 \le C M^{2/3}\norm U_{H^3}^{1/3},\qquad
 \norm{\partial_{ii}U}_3
 \le C M^{1/3}\norm U_{H^3}^{2/3}.
\]
Both bounds hold in every dimension. They follow by applying
the corresponding estimates on coordinate lines and integrating.
Consequently,
\begin{equation}\label{eq:weighted-third-chain}
 I_{i,\varepsilon}
 \le C_\mu M^{2\mu-4}
       \norm{\partial_iU\,\partial_{ii}U}_2^2
 \le C_\mu M^{2\mu-2}\norm U_{H^3}^2.
\end{equation}
This estimate is uniform as $\varepsilon\downarrow0$. In
particular, it justifies the singular weighted integral without
assuming its integrability in advance.

To apply the estimate to $f$, mollify $f$ near the range of $U$,
obtaining $f_\varepsilon$, and replace it by
$f_\varepsilon-f_\varepsilon(0)$. The hypotheses imply that
$f''$ is locally absolutely continuous across zero and that its
distributional derivative is the locally integrable function
$f'''$. Elementary convolution estimates give
\[
 |f_\varepsilon'''(\xi)|
 \le C_\mu(\xi^2+\varepsilon^2)^{(\mu-3)/2},\qquad
 |f_\varepsilon^{(j)}(\xi)|
 \le C_\mu(|\xi|+\varepsilon)^{\mu-j},\quad j=1,2,
\]
on the required range. For $0<\varepsilon<M$, the pure third
derivative identity
\[
 \partial_i^3(f_\varepsilon(U))
 =f_\varepsilon'''(U)(\partial_iU)^3
 +3f_\varepsilon''(U)\partial_iU\,\partial_{ii}U
 +f_\varepsilon'(U)\partial_i^3U
\]
is controlled by \eqref{eq:weighted-third-chain} and the preceding
product inequalities. Lower derivatives are treated by the usual
first- and second-order chain rules. Pure third derivatives,
together with the $L^2$ norm, control the full $H^3(\RR^n)$ norm
by Fourier transformation. We obtain
\[
 \norm{f_\varepsilon(U)-f_\varepsilon(0)}_{H^3}
 \le C_\mu M^{\mu-1}\norm U_{H^3}.
\]
Weak compactness in $H^3$ and convergence in $L^2$ prove the
corresponding bound for $f(U)$. If $\mu\ge3$, the same conclusion
follows directly from
$|f'''(U)|\le C M^{\mu-3}$ and
$\norm{\partial_iU}_6^3\le C M^2\norm U_{H^3}$;
mollification justifies the chain rule if necessary.

Approximation extends the estimate to compactly supported real
$U\in H^3(\RR^n)\cap L^\infty(\RR^n)$.
A smooth bounded domain admits an extension operator bounded
both on $H^3$ and on $L^\infty$ \cite[Chapter~VI]{Stein}, with compact support
after multiplication by a fixed cutoff. Apply the whole-space estimate to such an extension
of $u$. Choosing $r_*$ small enough keeps its range inside
$(-r_0,r_0)$ and proves \eqref{eq:real-composition-tame}.
\end{proof}

\subsection{Optimal decay from three additional spatial derivatives}

For $1\le n\le7$, define the sufficient threshold for the local
vanishing condition
\begin{equation}\label{eq:spatial-order-threshold}
 \mu_3(n)=
 \begin{cases}
 5/2,&1\le n\le4,\\
 3,&n=5,\\
 4,&n=6,\\
 7,&n=7.
 \end{cases}
\end{equation}
These are sufficient thresholds for the argument in the fixed space
$\mathcal X_3$.

\begin{lemma}[The nonlinear source in $\mathcal X_3$]
\label{lem:spatial-source-estimate}
Assume $1\le n\le7$, the source hypotheses of
Lemma~\ref{lem:real-power-composition}, and $\mu>\mu_3(n)$.
There is a number $\ell$ satisfying
\begin{equation}\label{eq:spatial-interpolation-index}
 1<\ell<4,\qquad \ell>n/2,\qquad
 (4-\ell)(\mu-1)>3.
\end{equation}
For $n\le3$, one may take $\ell=2$.
With
\begin{equation}\label{eq:spatial-force-exponents}
 \nu_\mu=\frac{(4-\ell)(\mu-1)+1}{3},\qquad
 \eta_\mu=\frac{3\nu_\mu}{4}>1,
\end{equation}
the source $\mathcal N(Y)=(0,f(u),0,0,0)$ satisfies
\begin{equation}\label{eq:spatial-source-tame}
 \norm{\mathcal N(Y)}_{\mathcal X_3}
 \le C\norm Y_E^{\nu_\mu}\norm Y_{\mathcal X_3}^{\mu-\nu_\mu}
\end{equation}
for $Y\in\mathcal X_3$ with sufficiently small $L^\infty$ norm
of its first component.
\end{lemma}
\begin{proof}
For $n\le3$, the choice $\ell=2$ is allowed by $\mu>5/2$.
For $4\le n\le7$, an index
\[
 n/2<\ell<4-\frac3{\mu-1}
\]
exists precisely under the inequalities stated in
\eqref{eq:spatial-order-threshold}.

We first verify the boundary condition in the source space.
On a smooth domain, the Dirichlet power-domain characterization
\cite[Theorem~2.2]{Grubb} gives
\[
 \DD(A^2)=\{u\in H^4:u=\Delta u=0\text{ on }\partial\Omega\},
\quad
 \DD(A^{3/2})=\{g\in H^3:g=\Delta g=0\text{ on }\partial\Omega\},
\]
with equivalent Sobolev and power norms. If $u$ is a finite
Dirichlet spectral sum, it is smooth, and $f\in C^2$ gives
\[
 f(u)=0,\qquad
 \Delta f(u)=f'(u)\Delta u+f''(u)|\nabla u|^2=0
 \quad\text{on }\partial\Omega.
\]
For a general $u\in\DD(A^2)$, its finite spectral sums converge
in $H^4$ and, since $n\le7$, uniformly. Lemma~\ref{lem:real-power-composition}
bounds their compositions in $H^3$, and local Lipschitz continuity
of $f$ identifies the weak $H^3$ limit as $f(u)$ by strong $L^2$
convergence. The continuous trace maps on $H^3$ preserve both
zero traces in this weak limit. Thus $f(u)\in\DD(A^{3/2})$.
This argument does not require $H^4\hookrightarrow C^2$ in the
higher dimensions.

Lemma~\ref{lem:real-power-composition} and interpolation now give
\[
 \begin{split}
 \norm{\mathcal N(Y)}_{\mathcal X_3}
 &\le C\norm u_\infty^{\mu-1}\norm u_{H^3},\\
 \norm u_\infty
 &\le C\norm u_{H^1}^{(4-\ell)/3}
             \norm u_{H^4}^{(\ell-1)/3},\\
 \norm u_{H^3}
 &\le C\norm u_{H^1}^{1/3}\norm u_{H^4}^{2/3}.
 \end{split}
\]
Here $H^\ell\hookrightarrow L^\infty$ follows from
$\ell>n/2$. Since $\norm u_{H^1}\le C\norm Y_E$ and
$\norm u_{H^4}\le C\norm Y_{\mathcal X_3}$, these inequalities
prove \eqref{eq:spatial-source-tame}.
\end{proof}

\begin{theorem}[Optimal small-data decay in $\mathcal X_3$]
\label{thm:power-log-decay}
Let $\Omega\subset\RR^n$ be a bounded $C^\infty$ domain,
$1\le n\le7$. Assume the source hypotheses of
Lemma~\ref{lem:real-power-composition}, and $\mu>\mu_3(n)$.
There exist $\varepsilon_*>0$ and $C>0$ such that every real
$Y_0\in\mathcal X_3$ with
$\norm{Y_0}_{\mathcal X_3}\le\varepsilon_*$ generates a unique
global regular solution of \eqref{eq:spatial-nonlinear-system},
\[
 Y\in C([0,\infty);\mathcal X_3)\cap C^1([0,\infty);\HH),
\]
and
\begin{equation}\label{eq:spatial-nonlinear-decay}
 \sup_{t\ge0}\norm{Y(t)}_{\mathcal X_3}
 +\sup_{t\ge0}(1+t)^{3/4}\norm{Y(t)}_E
 \le C\norm{Y_0}_{\mathcal X_3}.
\end{equation}
Its first component remains in the interval where the source hypotheses
are imposed. If $F(\xi)=\int_0^\xi f(\zeta)\,d\zeta$, its nonlinear energy
\[
 \EE_f(Y)=\tfrac12\norm Y_E^2-\int_\Omega F(u)\,dx
\]
satisfies
\begin{equation}\label{eq:spatial-nonlinear-energy}
 \begin{split}
 \EE_f(Y(t))+D\int_0^t\norm{A^{1/2}z(s)}_2^2\,ds&=\EE_f(Y_0),\\
 0\le \EE_f(Y(t))&\le
 C(1+t)^{-3/2}\norm{Y_0}_{\mathcal X_3}^2.
 \end{split}
\end{equation}
The exponents $3/4$ and $3/2$ are optimal for uniform estimates
proportional to the corresponding initial norms on a fixed small
$\mathcal X_3$ ball.
\end{theorem}
\begin{proof}
Choose a $C^2$ cutoff extension $\widetilde f$ agreeing with $f$
on a fixed neighborhood of zero and with bounded derivative on
$\RR$. Then
$\widetilde{\mathcal N}(Y)=(0,\widetilde f(u),0,0,0)$ is globally
Lipschitz from $\HH$ to $\HH$. The cutoff equation has a unique
global mild energy solution. It will follow from the estimates
below that the cutoff is never active for the initial states under
consideration.

Let $\Pi_N$ denote componentwise projection onto the first $N$
Dirichlet eigenfunctions. Its norm is at most one on both
$\HH$ and $\mathcal X_3$, and it commutes with $T(t)$.
Consider the finite-dimensional equations
\[
 (Y_N)_t=\GG Y_N+\Pi_N\widetilde{\mathcal N}(Y_N),
 \qquad Y_N(0)=\Pi_NY_0.
\]
They are globally defined. On any time interval on which their
first components stay in the amplitude range of
Lemma~\ref{lem:real-power-composition}, put
\[
 M_N(T)=\sup_{0\le t\le T}\norm{Y_N(t)}_{\mathcal X_3}
       +\sup_{0\le t\le T}(1+t)^{3/4}\norm{Y_N(t)}_E.
\]
Lemma~\ref{lem:spatial-source-estimate} implies
\begin{equation}\label{eq:spatial-force-time}
 \norm{\mathcal N(Y_N(t))}_{\mathcal X_3}
 \le C M_N(T)^\mu(1+t)^{-\eta_\mu},
 \qquad \eta_\mu>1.
\end{equation}
Duhamel's formula, contractivity on $\mathcal X_3$, and
\eqref{eq:spatial-decay} give
\[
 \begin{split}
 \norm{Y_N(t)}_{\mathcal X_3}
 &\le\norm{Y_0}_{\mathcal X_3}
 +CM_N(T)^\mu\int_0^t(1+s)^{-\eta_\mu}\,ds,\\
 \norm{Y_N(t)}_E
 &\le C(1+t)^{-3/4}\norm{Y_0}_{\mathcal X_3}
 +CM_N(T)^\mu
 \int_0^t(1+t-s)^{-3/4}(1+s)^{-\eta_\mu}\,ds.
 \end{split}
\]
Splitting the last integral at $t/2$ yields
\[
 \int_0^t(1+t-s)^{-3/4}(1+s)^{-\eta_\mu}\,ds
 \le C(1+t)^{-3/4},
\]
because $\eta_\mu>1$. Therefore
\begin{equation}\label{eq:spatial-bootstrap}
 M_N(T)\le C_0\norm{Y_0}_{\mathcal X_3}+C_1M_N(T)^\mu,
\end{equation}
with constants independent of $N$ and $T$.
The embedding $H^4(\Omega)\hookrightarrow L^\infty(\Omega)$
holds for $n\le7$. Choose $\varepsilon_*$ so that the usual
continuity bootstrap in \eqref{eq:spatial-bootstrap} gives
\[
 M_N(T)\le2C_0\norm{Y_0}_{\mathcal X_3},
\]
and so that this bound keeps $u_N$ strictly inside both the
source interval and the amplitude range of
Lemma~\ref{lem:real-power-composition}. The estimate consequently
holds for every $T$ and every $N$, and the cutoff is inactive.

Let $Y$ be the energy solution of the cutoff problem. Global
Lipschitz continuity of $\widetilde{\mathcal N}$, its Duhamel
formula, and Gronwall's inequality show that
\[
 Y_N\longrightarrow Y\quad\text{in }C([0,T];\HH)
 \qquad\text{for every }T<\infty.
\]
For clarity, the projection error tends uniformly to zero on
$[0,T]$, because
$\widetilde{\mathcal N}(Y([0,T]))$ is compact in $\HH$ and
$\Pi_N\to I$ strongly there. The uniform $\mathcal X_3$ bound
passes to the limit weakly at every time, while the energy decay
passes by strong convergence. Hence $Y(t)\in\mathcal X_3$,
the cutoff is inactive, and
\eqref{eq:spatial-nonlinear-decay} holds.

It remains to justify strong continuity in the higher space;
local Lipschitz continuity of
$\mathcal N:\mathcal X_3\to\mathcal X_3$ is not assumed.
The source $\mathcal N(Y(t))$ belongs pointwise to
$\mathcal X_3$ and satisfies the analogue of
\eqref{eq:spatial-force-time}. Its finite spectral projections
are continuous with values in $\mathcal X_3$, by energy
continuity of $Y$ and local Lipschitz continuity of the cutoff
source on $\HH$. They converge pointwise in $\mathcal X_3$ to
$\mathcal N(Y(t))$. Thus the latter is strongly measurable and
locally Bochner integrable in $\mathcal X_3$. Its Duhamel formula
in that space proves
$Y\in C([0,\infty);\mathcal X_3)$.
Since $\mathcal X_3\hookrightarrow\DD(\GG)$ and
$\GG:\mathcal X_3\to\HH$ is bounded, the equation gives
$Y\in C^1([0,\infty);\HH)$.
Uniqueness in the stated regular class follows in the energy
norm from local Lipschitz continuity of $f$ and Gronwall's
inequality while both solutions remain in the prescribed source
interval. A competing regular solution with the same initial
state cannot first leave this interval: it agrees up to that
time with the constructed solution, whose amplitude remains
strictly inside it. Thus uniqueness holds on every finite time
interval.

The source hypothesis implies $|F(\xi)|\le C|\xi|^{\mu+1}$ near zero.
Consequently,
\[
 \left|\int_\Omega F(u)\,dx\right|
 \le C\norm u_\infty^{\mu-1}\norm u_2^2
 \le C\norm u_\infty^{\mu-1}\norm Y_E^2.
\]
After a further reduction of $\varepsilon_*$,
\begin{equation}\label{eq:spatial-energy-equivalence}
 \tfrac14\norm{Y(t)}_E^2
 \le \EE_f(Y(t))\le\tfrac34\norm{Y(t)}_E^2.
\end{equation}
The regularity already proved justifies the energy identity
$\EE_f'(t)=-D\norm{A^{1/2}z(t)}_2^2$.
This proves \eqref{eq:spatial-nonlinear-energy}.

Finally, the same Duhamel estimates show that
\begin{equation}\label{eq:spatial-nonlinear-remainder}
 \norm{Y(t;Y_0)-T(t)Y_0}_E
 \le C\norm{Y_0}_{\mathcal X_3}^\mu(1+t)^{-3/4}.
\end{equation}
For fixed $Z\in\mathcal X_3$, substitute $Y_0=\delta Z$
and divide by $\delta$. Since $\mu>1$, letting $\delta\to0$
gives $\delta^{-1}Y(t;\delta Z)\to T(t)Z$ in $\HH$.
Any uniform state estimate on a small $\mathcal X_3$ ball with
an exponent greater than $3/4$, or a common additional factor
tending to zero, would therefore imply the same linear estimate,
contradicting Lemma~\ref{lem:spatial-decay}. By
\eqref{eq:spatial-energy-equivalence}, a faster uniform energy
rate would also give a faster state rate. This proves optimality.
\end{proof}

\begin{corollary}[Power and logarithmic examples]
\label{cor:power-log-spatial}
For $1\le n\le7$, Theorem~\ref{thm:power-log-decay} applies to
\begin{enumerate}
\item $f(\xi)=\kappa_0|\xi|^q \xi$, $\kappa_0>0$, provided
      $q>\mu_3(n)-1$;
\item $f(\xi)=|\xi|^{\gamma-2}\xi\log|\xi|$, $f(0)=0$, provided
      $\gamma>1+\mu_3(n)$.
\end{enumerate}
In particular, the sufficient ranges in the fixed space
$\mathcal X_3$ are
\[
\begin{array}{c|c|c}
 \text{dimension}&\text{power exponent }q&\text{logarithmic exponent }\gamma\\ \hline
 1\le n\le4&q>3/2&\gamma>7/2\\
 n=5&q>2&\gamma>4\\
 n=6&q>3&\gamma>5\\
 n=7&q>6&\gamma>8.
\end{array}
\]
\end{corollary}
\begin{proof}
For the power source, take $\mu=q+1$ in
\eqref{eq:spatial-source-assumption}. The source is $C^2$
because $q+1>5/2$, and the third derivative estimate is imposed
only away from zero when the source is not smooth there.
For the logarithmic source choose
\[
 \mu_3(n)<\mu<\gamma-1.
\]
For $j=0,1,2,3$ and $\xi\ne0$, its derivatives near zero satisfy
\[
 |f^{(j)}(\xi)|
 \le C|\xi|^{\gamma-1-j}(1+|\log|\xi||)
 \le C_\mu|\xi|^{\mu-j}.
\]
The derivatives of order at most two extend continuously through
zero and vanish there. Thus all source hypotheses hold.
\end{proof}

\begin{remark}[Scope and relation to the companion source ranges]
\label{rem:power-log-source-scope}
The preceding theorems give uniform decay from small regular initial
data, with the boundary compatibility encoded in their data spaces.
The $\mathcal X_3$ result also requires a smooth domain. These data
conditions supply the higher-order bound and the integrable source
estimate used in the proof.

In dimensions one and two, the power and logarithmic examples
belong to the source classes of \cite{HaGeneral}. In dimension
three, their overlap with its strictly $L^2$-subcritical ranges is
\[
 \frac32<q<2\quad\text{for }f(\xi)=\kappa_0|\xi|^q\xi,
 \qquad
 \frac72<\gamma<4\quad\text{for the logarithmic source}.
\]
For a pure power, the local vanishing exponent is $\mu=q+1$.
The cubic case $q=2$ is also covered here, although it is the
endpoint excluded by the companion paper's strict source range.
For $n\ge4$, the ranges in Corollary~\ref{cor:power-log-spatial}
lie outside the companion paper's strict $L^2$-subcritical power
and logarithmic ranges. They are small-data results in the
higher-regularity space $\mathcal X_3$; the amplitude bound ensures
that the cutoff construction solves the original equation.

The restriction $n\le7$ comes from $H^4\hookrightarrow L^\infty$.
For $n\ge8$, a higher spatial scale restores this embedding, as
in Theorem~\ref{thm:highdim-decay}. The source thresholds above are sufficient
conditions arising from the proofs. Decay for the remaining lower-order
sources, and unconditional strong convergence for every stable-well
trajectory, are not resolved by the present results.
\end{remark}

\subsection{High-order sources in arbitrary dimension}
\label{subsec:highdim-extension}

The next result covers every space dimension by requiring a higher
order of vanishing at zero. In dimensions $n\ge3$, these sources
include pure-power and logarithmic examples outside the strictly
$L^2$-subcritical ranges of the companion paper. The initial data also have stronger spatial
regularity.

\begin{theorem}\label{thm:highdim-decay}
Let $\Omega\subset\mathbb R^n$ be a bounded $C^\infty$ domain,
$n\ge1$, and let $m$ be an integer satisfying
\[
m\ge4,\qquad m\ge n/2.
\]
Suppose $f\in C^1(\mathbb R)$ and there exist $\mu>m$ and
$R>0$ such that $f$ is $C^m$ on an open neighborhood of
$[-2R,2R]$ and
\begin{equation}\label{eq:highdim-source}
|f^{(j)}(\xi)|\le C|\xi|^{\mu-j},\qquad
|\xi|\le2R,\quad 0\le j\le m.
\end{equation}
There exist $\varepsilon_*>0$ and $C_*>0$ such that every
$Y_0\in\mathcal X_m$ with
$\|Y_0\|_{\mathcal X_m}\le\varepsilon_*$ generates a unique
global solution
\[
Y\in C([0,\infty);\mathcal X_m)
\cap C^1([0,\infty);\mathcal H)
\]
of \eqref{eq:spatial-nonlinear-system}, satisfying
\begin{equation}\label{eq:highdim-decay}
\sup_{t\ge0}\|Y(t)\|_{\mathcal X_m}
+\sup_{t\ge0}(1+t)^{m/4}\|Y(t)\|_E
\le C_*\|Y_0\|_{\mathcal X_m}.
\end{equation}
Its nonlinear energy is equivalent to $\|Y(t)\|_E^2$ along the
solution and obeys
\[
\EE_f(Y(t))\le C(1+t)^{-m/2}\|Y_0\|_{\mathcal X_m}^2.
\]
The exponents $m/4$ and $m/2$ are optimal in these uniform
estimates on a fixed small $\mathcal X_m$ ball.
\end{theorem}
\begin{proof}
Choose a smooth cutoff equal to one on $[-R,R]$ and supported
in $(-2R,2R)$. The cutoff source, denoted by $\widetilde f$,
is globally $C^m$, has bounded first derivative, and satisfies
\eqref{eq:highdim-source} on all of $\mathbb R$ with a possibly
different constant. Its source map is globally Lipschitz on
$\mathcal H$.

For integer $m$, the Dirichlet power domain is
\cite[Theorem~2.2]{Grubb}
\[
\mathcal D(A^{m/2})
=\{g\in H^m(\Omega):\Delta^j g|_{\partial\Omega}=0
\text{ whenever }2j<m-\tfrac12\},
\]
with equivalent norms. Since $\mu>m$, all derivatives of
$\widetilde f$ through order $m$ vanish at zero. For a smooth
finite spectral sum $u$, each term of
$\Delta^j\widetilde f(u)|_{\partial\Omega}$ contains a factor
$\widetilde f^{(i)}(0)$ with $i\le2j\le m-1$.
Thus all required traces vanish. The integer-order tame
composition estimate (see \cite[Chapter~13]{Taylor}) gives
\begin{equation}\label{eq:highdim-tame}
\|A^{m/2}\widetilde f(u)\|_2
\le C\|\widetilde f(u)\|_{H^m}
\le C\|u\|_\infty^{\mu-1}\|u\|_{H^m}.
\end{equation}
The Sobolev estimate follows from the ordinary chain rule,
the derivative bounds in \eqref{eq:highdim-source}, and
Gagliardo--Nirenberg interpolation for each product of
derivatives of $u$. Spectral approximation in $H^{m+1}$, and
hence in $L^\infty$, extends the composition and its boundary
conditions to $u\in\mathcal D(A^{(m+1)/2})$.

Because $m\ge n/2$ and $\mu>m$, one can choose $\ell$ with
\[
\max\{1,n/2\}<\ell<m+1,
\qquad (m+1-\ell)(\mu-1)>m-1.
\]
Sobolev embedding and interpolation imply
\[
\begin{split}
\|u\|_\infty
&\le C\|Y\|_E^{(m+1-\ell)/m}
\|Y\|_{\mathcal X_m}^{(\ell-1)/m},\\
\|u\|_{H^m}
&\le C\|Y\|_E^{1/m}
\|Y\|_{\mathcal X_m}^{(m-1)/m}.
\end{split}
\]
Put $\beta=m/4\ge1$ and
\[
\eta=\frac{(m+1-\ell)(\mu-1)+1}{4}>\beta.
\]
For the weighted norm
\[
M(t)=\sup_{0\le s\le t}
\bigl(\|Y(s)\|_{\mathcal X_m}
+(1+s)^\beta\|Y(s)\|_E\bigr),
\]
estimate \eqref{eq:highdim-tame} gives
\[
\|\widetilde{\mathcal N}(Y(s))\|_{\mathcal X_m}
\le CM(t)^\mu(1+s)^{-\eta}.
\]
In particular the source is time integrable in $\mathcal X_m$.
The elementary convolution bound
\[
\int_0^t(1+t-s)^{-\beta}(1+s)^{-\eta}\,ds
\le C(1+t)^{-\beta},\qquad \eta>\beta\ge1,
\]
follows by splitting at $t/2$; when $\beta=1$, the second half
has the harmless bound
$C(1+t)^{-\eta}\log(2+t)$.

The spectral Galerkin argument used above now applies with
$s=m$. The projections commute with the linear semigroup and
are contractions on $\mathcal X_m$. Duhamel's formula gives,
uniformly in the Galerkin dimension,
\[
M(t)\le C_0\|Y_0\|_{\mathcal X_m}+C_1M(t)^\mu.
\]
Smallness closes this inequality globally. The globally
Lipschitz cutoff source gives convergence in
$C([0,T];\mathcal H)$ for every finite $T$, while weak lower
semicontinuity gives the uniform $\mathcal X_m$ bound. As in
the preceding Galerkin proof, the spectral coordinates make
the source strongly measurable in the separable space
$\mathcal X_m$; its time-integrable bound turns Duhamel's
formula into a Bochner integral there. Thus the limit belongs
to $C([0,\infty);\mathcal X_m)$. No Lipschitz composition
estimate in $\mathcal X_m$ is needed. Since
$\mathcal X_m\subset\mathcal D(\mathcal G)$, the solution is
classical in $\mathcal H$.

The uniform $H^{m+1}$ bound keeps $\|u(t)\|_\infty\le R$
after decreasing $\varepsilon_*$. Consequently this solution
solves the original equation without cutoff. Uniqueness in
$C([0,T];\mathcal X_m)$ follows from the pointwise Lipschitz
bound for $f$ on the bounded ranges of the two displacements
and Gronwall's inequality in $\mathcal H$.
Writing $F(\xi)=\int_0^\xi f(\zeta)\,d\zeta$, we have
\[
\left|\int_\Omega F(u)\,dx\right|
\le C\|u\|_\infty^{\mu-1}\|u\|_2^2
\le C\varepsilon_*^{\mu-1}\|Y\|_E^2.
\]
This proves energy equivalence and decay even when the source
is energy-supercritical. The regularity above also justifies
$\EE_f'(t)=-D\|A^{1/2}z(t)\|_2^2$.
Finally, the same convolution estimate gives
\[
\|Y(t;Y_0)-T(t)Y_0\|_E
\le C\|Y_0\|_{\mathcal X_m}^{\mu}(1+t)^{-m/4}.
\]
Substituting $Y_0=\delta Z$, dividing by $\delta$, and letting
$\delta\to0$ identifies the tangent flow with $T(t)$.
Any better uniform nonlinear state rate would imply a better rate from
$\mathcal X_m$ to $\mathcal H$ for the linear semigroup,
contradicting Lemma~\ref{lem:spatial-decay}. Energy equivalence
then gives the asserted optimality of the energy exponent.
\end{proof}

\begin{corollary}[Pure powers and logarithmic powers]
For $f(\xi)=|\xi|^q \xi$, the theorem applies for every $q>m-1$.
For $f(\xi)=|\xi|^{\gamma-2}\xi\log|\xi|$, it applies for every
$\gamma>m+1$: choose $m<\mu<\gamma-1$ and use
$|\xi|^{\gamma-1-j}|\log|\xi||\le C|\xi|^{\mu-j}$ near zero for
$0\le j\le m$. Thus sufficiently high-order pure and
logarithmic powers satisfy the decay estimates above in
every dimension.
\end{corollary}

\section{Conditional decay under weighted variation}\label{sec:nonlinear-variation}

The abstract transfer result below retains the Hilbert-space setting of Section~\ref{sec:abstract}. Its PDE application is on a bounded $C^2$ domain in any dimension $n\ge1$. In the nonlinear application, the weighted-variation assumption is an additional trajectory condition; it is
not required in any of the preceding small-data theorems.

\begin{proposition}[Weighted variation of an inhomogeneous term]\label{prop:weighted-force}
Let $Y_0\in\DD(\GG)$, let
$g\in AC_{\mathrm{loc}}([0,\infty);\HH)$, and suppose that
\begin{equation}\label{eq:weighted-force-assumption}
 \lim_{t\to\infty}\norm{g(t)}_E=0,
 \qquad
 V_{1/4}(g):=\int_0^\infty(1+s)^{1/4}\norm{g'(s)}_E\,ds<\infty.
\end{equation}
The mild solution of
\[
 Y_t=\GG Y+g(t),\qquad Y(0)=Y_0,
\]
satisfies
\begin{equation}\label{eq:weighted-force-rate}
 \sup_{t\ge0}(1+t)^{1/4}\norm{Y(t)}_E
 \le C\bigl(\norm{\GG Y_0+g(0)}_E+V_{1/4}(g)\bigr).
\end{equation}
In fact,
\begin{equation}\label{eq:weighted-force-little-o}
 \norm{Y(t)}_E=o(t^{-1/4})\qquad(t\to\infty).
\end{equation}
The constant $C$ depends only on the linear system.
\end{proposition}

\begin{proof}
Since $0\in\rho(\GG)$, integration by parts in the variation-of-constants
formula gives
\begin{equation}\label{eq:forced-ibp}
 \begin{split}
 Y(t)={}&T(t)\bigl[Y_0+\GG^{-1}g(0)\bigr]-\GG^{-1}g(t)\\
 &+\int_0^tT(t-s)\GG^{-1}g'(s)\,ds.
 \end{split}
\end{equation}
To justify this identity under the stated assumptions, differentiate
$T(t-s)\GG^{-1}g(s)$ for almost every $s\in[0,t]$.
The function $\GG^{-1}g$ is locally absolutely continuous with values
in $\DD(\GG)$ equipped with its graph norm, and its derivative is
$\GG^{-1}g'$. Thus the product rule and Bochner integration apply.

The first assumption in \eqref{eq:weighted-force-assumption} and the
integrability of $g'$ imply
\[
 g(t)=-\int_t^\infty g'(s)\,ds,
 \qquad
 (1+t)^{1/4}\norm{g(t)}_E
 \le\int_t^\infty(1+s)^{1/4}\norm{g'(s)}_E\,ds.
\]
Moreover, for $0\le s\le t$,
\[
 1+t\le(1+t-s)(1+s).
\]
Applying the linear estimate
$\norm{T(r)\GG^{-1}}_E\le C(1+r)^{-1/4}$ to
\eqref{eq:forced-ibp} proves \eqref{eq:weighted-force-rate}.

For the little-$o$ assertion, the initial vector
$Y_0+\GG^{-1}g(0)$ belongs to $\DD(\GG)$, so its orbit satisfies
the fixed-datum conclusion of Proposition~\ref{prop:orbits}.
The weighted tail estimate above implies
$(1+t)^{1/4}\norm{\GG^{-1}g(t)}_E\to0$.
For almost every fixed $s$, the same fixed-datum conclusion gives
\[
 (1+t)^{1/4}T(t-s)\GG^{-1}g'(s)\longrightarrow0
 \quad\hbox{in }\HH.
\]
After extending this expression by zero when $s>t$, its norm is
bounded by $C(1+s)^{1/4}\norm{g'(s)}_E$.
Dominated convergence in the last term of \eqref{eq:forced-ibp}
therefore proves \eqref{eq:weighted-force-little-o}.
\end{proof}

We now apply Proposition~\ref{prop:weighted-force} to the general
focusing-source class in
\cite{HaGeneral}. Recall that its derivative-growth hypothesis is
\begin{equation}\label{eq:general-source-growth-variation}
 \begin{gathered}
 f\in C^1(\RR),\qquad f(0)=f'(0)=0,\qquad
 |f'(\xi)|\le C_f(1+|\xi|^\varrho),\\
 0<\varrho<\varrho_*,\qquad
 \varrho_*=\begin{cases}\infty,&n\le2,\\2/(n-2),&n\ge3.\end{cases}
 \end{gathered}
\end{equation}
The stable-set theory also assumes the nonnegativity of
$H_\theta(\xi)=\theta^{-1}\xi f(\xi)-F(\xi)$ for some $\theta>2$, where
$F(\xi)=\int_0^\xi f(\zeta)\,d\zeta$, and the existence of
$\phi\in H_0^1(\Omega)$ with $\int_\Omega F(\phi)\,dx>0$.
These are hypotheses (F1), (F2), and (F4) in that reference.
The additional radial monotonicity hypothesis (F3), used in its
unstable-side results, is not needed for the present application.

\begin{corollary}[Conditional decay of general stable trajectories]
\label{cor:general-weighted-variation}
Assume \eqref{eq:parameters} and hypotheses (F1), (F2), and (F4)
of \cite{HaGeneral}. Let $Y$ be the global regular stable solution
provided by Theorem~5.13 of that reference, with initial energy
strictly below the coupled well depth. In particular,
\[
 Y_0\in
 \bigl(H^2(\Omega)\cap H_0^1(\Omega)\bigr)
 \times H_0^1(\Omega)
 \times\bigl(H^2(\Omega)\cap H_0^1(\Omega)\bigr)
 \times H_0^1(\Omega)
 \times\bigl(H^2(\Omega)\cap H_0^1(\Omega)\bigr).
\]
If
\begin{equation}\label{eq:weighted-nonlinear-variation}
 \int_0^\infty(1+s)^{1/4}
 \norm{f'(u(s))u_t(s)}_{L^2(\Omega)}\,ds<\infty,
\end{equation}
then
\[
 \norm{Y(t)}_E=o(t^{-1/4}),\qquad
 \EE_f(t)=o(t^{-1/2}),
 \qquad
 \EE_f(t)=\tfrac12\norm{Y(t)}_E^2-\int_\Omega F(u(t))\,dx.
\]
The corresponding big-$O$ state estimate is
\eqref{eq:weighted-force-rate} with
$g(t)=(0,f(u(t)),0,0,0)$.
\end{corollary}

\begin{proof}
The regular phase space displayed above is contained in
$\DD(\GG)$. The regular chain rule in Appendix~A of
\cite{HaGeneral} gives
\[
 g\in C^1([0,\infty);\HH),\qquad
 g'(t)=(0,f'(u(t))u_t(t),0,0,0).
\]
The unconditional weak-stability result in Theorem~5.4 of that
reference yields $u(t)\to0$ in every strictly subcritical Lebesgue
space. In particular,
$u(t)\to0$ in $L^{2(\varrho+1)}(\Omega)$ because
$2(\varrho+1)<2^*$ when $n\ge3$; for $n\le2$, the exponent is finite.
Since \eqref{eq:general-source-growth-variation} implies
$|f(\xi)|\le C(|\xi|+|\xi|^{\varrho+1})$, it follows that
$g(t)\to0$ in $\HH$.
Proposition~\ref{prop:weighted-force} therefore applies.

Finally, the same growth condition gives
$|F(\xi)|\le C(|\xi|^2+|\xi|^{\varrho+2})$.
Sobolev embedding, $\varrho+2<2^*$ for $n\ge3$, and the state estimate
show that
\[
 |\EE_f(t)|\le C\bigl(\norm{Y(t)}_E^2+
                       \norm{Y(t)}_E^{\varrho+2}\bigr)
 =o(t^{-1/2}).
\]
\end{proof}

\begin{remark}
The unweighted condition
$\int_0^\infty\norm{f'(u)u_t}_2\,dt<\infty$
in \cite{HaGeneral} is sufficient for strong convergence and a
uniform regular bound. It gives no specified algebraic decay of
the variation tail of the nonlinear source. Condition
\eqref{eq:weighted-nonlinear-variation} supplies this missing
quantitative information.
\end{remark}

\appendix
\section{An application to coupled waves with unequal speeds}\label{app:wave-wave}

For coupled waves with unequal propagation speeds, Oquendo and Raya \cite{OR} established the optimal semigroup decay rate $t^{-1/4}$. Here Proposition~\ref{prop:modal-principle} yields the exact resolvent and semigroup asymptotic constants in the specified energy norm. Retain $X$, $A$, and $V$ from Section~\ref{sec:abstract}, and consider
\begin{equation}\label{eq:ww-system}
u_{tt}+Au+\chi v=0,\qquad
v_{tt}+c_0^2Av+\chi u+d_0v_t=0,
\end{equation}
where $d_0>0$, $c_0>0$, $c_0\ne1$, and $\chi\in\RR$ satisfies $0<|\chi|<c_0\lambda_1$.
For $Z=(u,p_u,v,p_v)$ in $\HH_{\mathrm{ww}}=V\times X\times V\times X$, define
\begin{equation}\label{eq:ww-energy}
\norm{Z}_{E,\mathrm{ww}}^2
=\norm{A^{1/2}u}^2+\norm{p_u}^2
+c_0^2\norm{A^{1/2}v}^2+\norm{p_v}^2
+2\chi\Rea\ip{u}{v}.
\end{equation}
The condition on $\chi$ makes this norm equivalent to the product norm. Set
\[
\mathcal L_{\mathrm{ww}}(u,p_u,v,p_v)
=(p_u,-Au-\chi v,p_v,-c_0^2Av-\chi u-d_0p_v),
\]
with domain $\DD(A)\times V\times\DD(A)\times V$.

\begin{proposition}\label{prop:ww-example}
The operator $\mathcal L_{\mathrm{ww}}$ generates a contraction semigroup $\mathcal T_{\mathrm{ww}}(t)$ in the norm \eqref{eq:ww-energy}, with $i\RR\subset\rho(\mathcal L_{\mathrm{ww}})$. Setting $a_{\mathrm{ww}}=\chi^2d_0/[2(c_0^2-1)^2]$, one has
\begin{equation}\label{eq:ww-res}
\limsup_{|\omega|\to\infty}|\omega|^{-4}
\norm{(i\omega-\mathcal L_{\mathrm{ww}})^{-1}}_{E,\mathrm{ww}}
=a_{\mathrm{ww}}^{-1}.
\end{equation}
Moreover, for every integer $m\ge1$,
\begin{equation}\label{eq:ww-time}
\limsup_{t\to\infty}t^{m/4}
\norm{\mathcal T_{\mathrm{ww}}(t)\mathcal L_{\mathrm{ww}}^{-m}}_{E,\mathrm{ww}}
=\left(\frac{m}{4ea_{\mathrm{ww}}}\right)^{m/4}.
\end{equation}
If $\lambda_{k+1}/\lambda_k\to1$, the upper limit in \eqref{eq:ww-time} is a limit.
\end{proposition}
\begin{proof}
For the mode $r=\sqrt{\lambda_k}$, use coordinates $(ru_k,p_{u,k},c_0rv_k,p_{v,k})$, with metric
$G_{\mathrm{ww}}(r)=I+\chi(c_0r^2)^{-1}(E_{13}+E_{31})$.
$\mathcal L_{\mathrm{ww}}$ is the maximal direct sum of the modal matrices $B_{\mathrm{ww}}(r_k)$: the first and third components require $p_u,p_v\in V$, and the others require $u,v\in\DD(A)$. Each block is dissipative, since
$\frac{d}{dt}\norm{Z(t)}_{E,\mathrm{ww}}^2=-2d_0\norm{p_v(t)}^2$.
A purely imaginary eigenvalue forces $p_v=0$. At a nonzero imaginary eigenvalue the modal equations then give $v=u=p_u=0$; at zero the remaining system is invertible because $c_0^2r^4-\chi^2>0$. Thus all blocks are strictly stable.

With $\varepsilon=r^{-1}$ and $\widehat{s}=s/r$, the characteristic equation becomes
\[
(\widehat{s}^2+1)(\widehat{s}^2+c_0^2+d_0\varepsilon\widehat{s})-\chi^2\varepsilon^4=0.
\]
At $\varepsilon=0$, the roots $\pm i,\pm ic_0$ are simple. The analytic implicit function theorem and coefficient comparison give, for $\sigma=\pm1$,
\begin{align*}
s_{\mathrm s,\sigma}(r)
&=\sigma ir-\frac{\sigma i\chi^2}{2(c_0^2-1)}r^{-3}
-a_{\mathrm{ww}}r^{-4}+O(r^{-5}),\\
s_{\mathrm f,\sigma}(r)&=\sigma ic_0r-\frac{d_0}{2}+O(r^{-1}).
\end{align*}
Moreover, $G_{\mathrm{ww}}(r)\to I$ and
$r^{-1}B_{\mathrm{ww}}(r)\to\diag(J_2,c_0J_2)$,
where $J_2=\left(\begin{smallmatrix}0&1\\-1&0\end{smallmatrix}\right)$.
Energy normalization preserves this skew-Hermitian limit with simple spectrum, allowing diagonalizers with condition numbers $\kappa_k\to1$ in \eqref{eq:general-diagonal}.
Proposition~\ref{prop:modal-principle} with $\beta=4$, $c_\pm=\pm1$, and $L_0=a_{\mathrm{ww}}^{-1}$ gives the claims, with contractivity following from modal dissipativity.
\end{proof}

\section*{Acknowledgments}
This research was supported by Basic Science Research Program through the National Research Foundation of Korea (NRF) funded by the Ministry of Education (RS-2022-NR075641).


\begin{thebibliography}{99}

\bibitem{ACK}
F.~Alabau, P.~Cannarsa, and V.~Komornik,
Indirect internal stabilization of weakly coupled evolution equations,
\emph{J. Evol. Equ.} \textbf{2} (2002), 127--150.

\bibitem{AL}
F.~Alabau-Boussouira and M.~L\'eautaud,
Indirect stabilization of locally coupled wave-type systems,
\emph{ESAIM Control Optim. Calc. Var.} \textbf{18} (2012), 548--582.

\bibitem{BD}
C.~J.~K. Batty and T.~Duyckaerts,
Non-uniform stability for bounded semi-groups on Banach spaces,
\emph{J. Evol. Equ.} \textbf{8} (2008), 765--780.

\bibitem{BT}
A.~Borichev and Y.~Tomilov,
Optimal polynomial decay of functions and operator semigroups,
\emph{Math. Ann.} \textbf{347} (2010), 455--478.

\bibitem{BourdaudSharp}
G.~Bourdaud,
A sharpness result for powers of Besov functions,
\emph{J. Funct. Spaces Appl.} \textbf{2} (2004), 267--277.

\bibitem{Bourdaud}
G.~Bourdaud,
An introduction to composition operators in Sobolev spaces,
\emph{Eurasian Math. J.} \textbf{14} (2023), no.~1, 39--54.

\bibitem{BucciEller}
F.~Bucci and M.~Eller,
The Cauchy--Dirichlet problem for the Moore--Gibson--Thompson equation,
\emph{C. R. Math. Acad. Sci. Paris} \textbf{359} (2021), 881--903.

\bibitem{Davies}
E.~B. Davies,
\emph{Spectral Theory and Differential Operators},
Cambridge Studies in Advanced Mathematics, vol.~42,
Cambridge University Press, Cambridge, 1995.

\bibitem{FQ}
J.~R. Fern\'andez and R.~Quintanilla,
On a mixture of an MGT viscous material and an elastic solid,
\emph{Acta Mech.} \textbf{233} (2022), 291--297.

\bibitem{Grubb}
G.~Grubb,
Regularity of spectral fractional Dirichlet and Neumann problems,
\emph{Math. Nachr.} \textbf{289} (2016), 831--844.

\bibitem{HaJDE}
T.~G. Ha,
Global well-posedness and conditional asymptotic stability for a coupled wave--MGT system with logarithmic nonlinearity,
\emph{J. Differential Equations} \textbf{485} (2026), 114748.

\bibitem{HaGeneral}
T.~G. Ha,
Stable and unstable potential-well dynamics for an indirectly damped wave--MGT system with general focusing sources,
preprint, 2026, version~1, \href{https://arxiv.org/abs/2609.13708v1}{arXiv:2609.13708v1}.

\bibitem{JL}
R.~Joly and C.~Laurent,
Decay of semilinear damped wave equations: cases without geometric control condition,
\emph{Ann. H. Lebesgue} \textbf{3} (2020), 1241--1289.

\bibitem{KLM}
B.~Kaltenbacher, I.~Lasiecka, and R.~Marchand,
Wellposedness and exponential decay rates for the Moore--Gibson--Thompson equation arising in high intensity ultrasound,
\emph{Control Cybernet.} \textbf{40} (2011), 971--988.

\bibitem{Kateb}
D.~Kateb,
On the boundedness of the mapping $f\mapsto |f|^\mu$, $\mu>1$, on Besov spaces,
\emph{Math. Nachr.} \textbf{248--249} (2003), 110--128.

\bibitem{MMT}
R.~Marchand, T.~McDevitt, and R.~Triggiani,
An abstract semigroup approach to the third-order Moore--Gibson--Thompson partial differential equation arising in high-intensity ultrasound: structural decomposition, spectral analysis, exponential stability,
\emph{Math. Methods Appl. Sci.} \textbf{35} (2012), 1896--1929.

\bibitem{Nirenberg}
L.~Nirenberg,
On elliptic partial differential equations,
\emph{Ann. Scuola Norm. Sup. Pisa Cl. Sci. (3)} \textbf{13} (1959), no.~2, 115--162.

\bibitem{OOB}
R.~L. Oliveira, H.~P. Oquendo, and C.~Buriol,
Exact decay rates for weakly coupled systems with indirect damping by fractional memory effects,
\emph{Math. Nachr.} \textbf{296} (2023), 4610--4633.

\bibitem{OR}
H.~P. Oquendo and R.~P. Raya,
Best rates of decay for coupled waves with different propagation speeds,
\emph{Z. Angew. Math. Phys.} \textbf{68} (2017), no.~4, Paper No.~77.

\bibitem{Pazy}
A.~Pazy,
\emph{Semigroups of Linear Operators and Applications to Partial Differential Equations},
Applied Mathematical Sciences, vol.~44,
Springer-Verlag, New York, 1983.

\bibitem{Stein}
E.~M. Stein,
\emph{Singular Integrals and Differentiability Properties of Functions},
Princeton Mathematical Series, vol.~30,
Princeton University Press, Princeton, NJ, 1970.

\bibitem{Taylor}
M.~E. Taylor,
\emph{Partial Differential Equations III: Nonlinear Equations},
2nd ed., Applied Mathematical Sciences, vol.~117,
Springer, New York, 2011.

\bibitem{XL}
T.-J. Xiao and J.~Liang,
Coupled second order semilinear evolution equations indirectly damped via memory effects,
\emph{J. Differential Equations} \textbf{254} (2013), 2128--2157.

\end{thebibliography}
\end{document}